\documentclass[11pt]{article}

\usepackage{amsmath, amsthm, amssymb}
\usepackage{mathrsfs}
\usepackage{graphics,amsopn,amstext,amsfonts,color}
\usepackage{bm}
\usepackage{fullpage}
\usepackage{setspace}
\usepackage{subfigure}
\usepackage{graphicx}
\usepackage{epstopdf}
\usepackage{natbib}
\usepackage{caption}
\usepackage{multirow}
\usepackage{booktabs}
\usepackage{paralist}
\usepackage{colortbl}
\usepackage{appendix}
\usepackage[bookmarksnumbered=true,bookmarksnumbered=true,
bookmarksopen=true]{hyperref}
\DeclareCaptionType{mycapequ}[][List of equations]
\bibpunct{(}{)}{;}{a}{,}{,}
\usepackage{setspace}
\usepackage{chngcntr}
\usepackage{enumerate}
\makeatletter
\renewcommand\paragraph{\@startsection{paragraph}{4}{\z@}%
{-3.25ex\@plus -1ex \@minus -.2ex}%
{0.0001pt \@plus .2ex}%
{\normalfont\normalsize\bfseries}}
\renewcommand\subparagraph{\@startsection{subparagraph}{5}{\z@}%
{-3.25ex\@plus -1ex \@minus -.2ex}%
{0.0001pt \@plus .2ex}%
{\normalfont\normalsize\bfseries}}
\counterwithin{paragraph}{subsubsection}
\counterwithin{subparagraph}{paragraph}
\makeatother
\theoremstyle{plain}

\newtheorem{corollary}{Corollary}
\newtheorem{lemma}{Lemma}

\newtheorem{proposition}{Proposition}
\newtheorem{property}{Property}

\theoremstyle{definition}
\newtheorem{definition}{Definition}
\theoremstyle{assumption}

\theoremstyle{remark}

\usepackage{array}
\newcolumntype{L}{>{\centering\arraybackslash}m{0.5cm}}
\usepackage{authblk}
\title{Optimal Layout of Transshipment Facilities on An Infinite Homogeneous Plane}

\author[1]{Weijun Xie }
\author[2]{Yanfeng Ouyang \thanks{Corresponding author, email: yfouyang@illinois.edu.}}
\affil[1]{School of Industrial \& Systems Engineering, Georgia Institute of Technology, Atlanta, GA 30332}
\affil[2]{Department of Civil and Environmental Engineering, University of Illinois, Urbana-Champaign, Urbana, IL 61801}

\date{}

\begin{document}

\maketitle

\begin{abstract}
This paper studies optimal spatial layout of transshipment facilities and the corresponding service regions on an infinite homogeneous plane $\Re^2$ that minimize the total cost for facility set-up, outbound delivery and inbound replenishment transportation. The problem has strong implications in the context of freight logistics and transit system design. This paper first focuses on a Euclidean plane and presents a new proof for the known Gersho's conjecture, which states that the optimal shape of each service region should be a regular hexagon if the inbound transportation cost is ignored. When inbound transportation cost becomes non-negligible, however, we show that a tight upper bound can be achieved by a type of elongated cyclic hexagons, while a cost lower bound based on relaxation and idealization is also obtained. The gap between the analytical upper and lower bounds is within 0.3\%. This paper then shows that a similar elongated non-cyclic hexagon shape is actually optimal for service regions on a rectilinear metric plane. Numerical experiments and sensitivity analyses are conducted to verify the analytical findings and to draw managerial insights.
\end{abstract}


\section{Introduction}\label{Intro}
Problems related to facility location (e.g., fixed-charge location problems) and routing (e.g., travel salesman problem, or TSP) impose two fundamental yet distinct challenges to logistics system design. Due to their intrinsic complexity, these problems are typically handled separately in the literature (e.g., see \citealp{daskin95} and \citealp{toth2001vehicle} for complete reviews). Relatively fewer studies looked at the integrated ``location-routing'' problem with or without inventory considerations (e.g., \citealp{Perl1985,Shen2007incorporating}). The location of distribution centers and the routing of outbound delivery vehicles are optimized simultaneously while inbound shipment (i.e., providing replenishment to these facilities) is assumed to be via direct visits (or more often, omitted from the model). Furthermore, most of these efforts focused on developing discrete mathematical programming models which can only numerically solve very limited-scale problem instances. In particular, little is known about the optimal facility lay-out, the suitable customer allocation, and the optimal vehicle tour in infinite homogeneous planes. Recently, \citet{Cachon11} proposed a new location-routing model in a homogeneous Euclidean plane which optimizes the spatial layout of facilities that serve distributed customers (i.e., similar to a median problem) and the routing of an inbound vehicle which visits these facilities (i.e., similar to a TSP problem). The model tried to minimize the total cost related to outbound customer access (i.e., direct shipment) and inbound replenishment transportation.

This problem has strong implications on practical logistics systems design in real-world contexts. For example, transshipment is often used in a timber harvesting system, where the lumbers are collected by trucks to local processing mills, and then shipped out by train. Train capacity is generally orders of magnitudes larger than that of local trucks, and it is often sufficient to assume infinite capacity of a train and design a single train track route for a large area of forest. Another example is commuter transit system design in low-demand areas~\citep*{nourbakhsh2012structured}, where a bus route collects passengers from a certain region (with sufficient capacity) at optimally located bus stops (where passengers gather).

\begin{figure}[htbp]
 \centering
 \subfigure[Inbound and outbound transportation]{
 \label{figure_notation} 
 \includegraphics[width=0.4\textwidth]{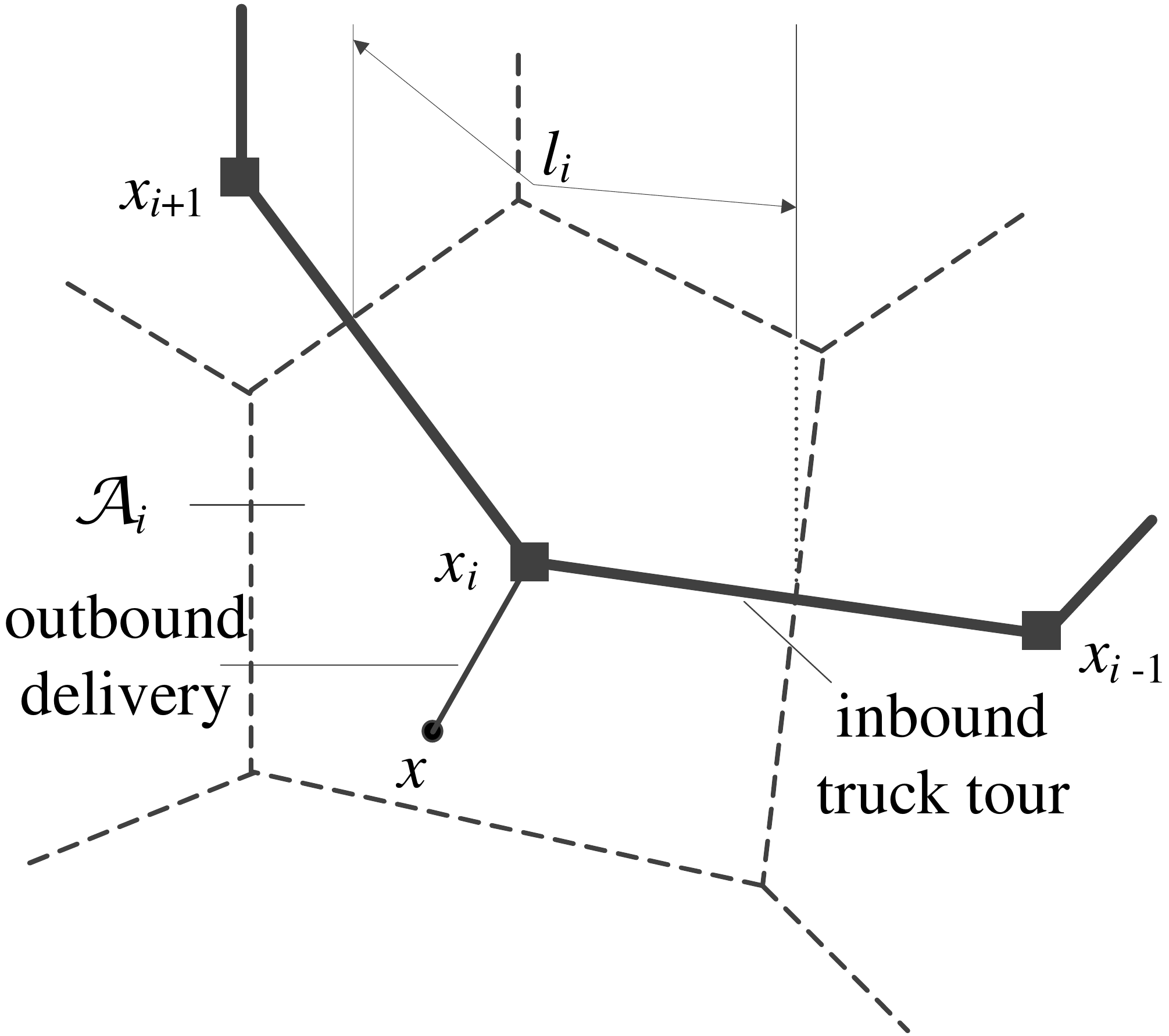}}
 \hspace{0.3in}
 \subfigure[Basic triangle and basic angle]{
 \label{Fig_Basic_Triangle}
 \includegraphics[width=0.4\textwidth]{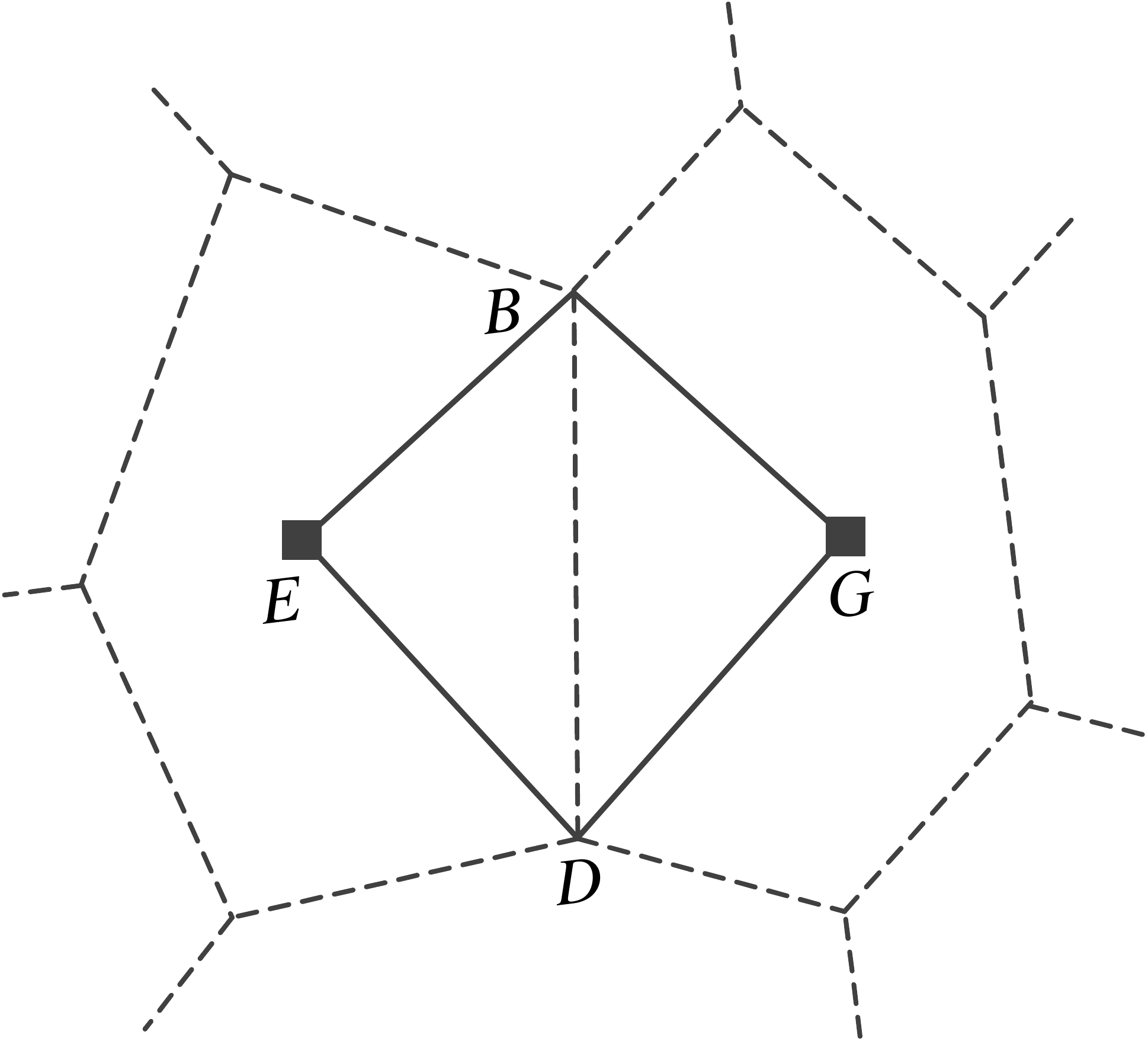}}
 \caption{Notations for a transshipment system}
 \label{Fig_Noatation} 
\end{figure}

Mathematically, this problem can be described as follows. We use a transshipment system to serve uniformly distributed customers (with demand density $\lambda$ per area-time) on an infinite homogeneous Euclidean plane $\Re^2$. Transshipment facilities can be constructed anywhere with a prorated set-up and operational cost $f$ per facility-time. All facilities receive replenishment from a central depot, which is co-located at one of the facilities, through an inbound truck with infinite capacity. This truck will supply all facilities along one tour, incurring average transportation cost of $C$ per distance\footnote{\noindent{This assumption makes the replenishment frequency irrelevant to the optimization problem.}}. Without losing generality, we assume that the transshipment facilities are indexed along the tour of the inbound truck; i.e., the inbound truck starts from facility 1 (the depot), visits facilities sequentially in set ${\cal{N}}=\left\{1,2,\cdots, N \right\}$ before returning to facility 1, where $N:=\left\vert {\cal{N}}\right\vert\rightarrow \infty$ is the total number of facilities. Facility $i \in {\cal{N}}$ is located at $x_i\in\Re^2$ to serve the customers in its service region ${\cal{A}}_{i}\in\Re^2$ through direct shipment, with a transportation cost of $c$ per demand-distance. For each customer at $x \in {\cal{A}}_{i}$, we use $\|x-x_i\|$ to denote the outbound travel distance. Moreover, for simplicity, the size of service region ${\cal{A}}_{i}$ is denoted by 
$A_i=\left\vert {\cal{A}}_{i}\right\vert$, and we assume the inbound truck travels a distance of $l_i$ within $A_i$.
Since all cost terms are relative, we further
define $\kappa=\frac{c\lambda}{f}$ (area-demand-distance/facility) and $r =\frac{C}{c\lambda}$ to denote the relative magnitudes of facility cost and inbound transportation cost as compared to the outbound cost.
Figure~\ref{figure_notation} illustrates these notations.

Some basic properties of this problem is readily available. For any given facility layout $\left\{x_i:i\in{\cal N}\right\}$, each customer should obviously choose the nearest facility for service and a tie may be broken arbitrarily. Thus, the set of service regions $\left\{{\cal A}_i:i\in{\cal N}\right\}$ must form a Voronoi diagram~\citep{Okabe1992,Du1999centroidal}, where
\begin{equation}
{\cal A}_i=\left\{x\in\Re^2:\|x-x_i\|\leq \|x-x_j\|,\forall j\neq i\right\}.\label{eq:A_i_def}
\end{equation}

Moreover, for an optimal TSP tour along $N$ facility locations in the Euclidean plane, it is easy to see that the optimal TSP tour has no crossover between any four facility locations; otherwise, a simple local perturbation can improve the solution. With this, the optimal solution of our problem must satisfy the following properties.
\begin{property}\label{P1} For any given facility locations $\left\{x_i:i \in {\cal{N}}\right\}$,
\begin{itemize}
\item the optimal service regions form a Voronoi diagram, i.e.,~\eqref{eq:A_i_def} holds and
 \begin{equation}
 \bigcup\limits_{i}{\cal A}_i=\Re^2; 
\label{eq_int_A_i_A_j}
 \end{equation}
\item each service region ${\cal A}_i$ is a convex polygon~\citep{Okabe1992} with the number of sides $n_i$;
\item for all $x\in {\cal A}_i\bigcap {\cal A}_{j}$, $\|x-x_i\|=\|x-x_{j}\|$, and therefore, for all $i \in {\cal{N}}$,\begin{equation}l_i=\frac{1}{2}\left(\|x_{i}-x_{i-1}\|+\|x_{i}-x_{i+1}\|\right);\label{eq:l_i_def}
\end{equation}
\item Line $x_ix_{i+1}$ is perpendicular to the interception line of ${\cal A}_i$ and ${\cal A}_{i+1}$.
\end{itemize}
\end{property}

In an infinite plane, a suitable objective is to find the optimal facility layout $\left\{x_i:i \in {\cal{N}}\right\}$, the service region partition $\left\{{\cal A}_i:i \in {\cal{N}}\right\}$, and the inbound truck tour that minimize the total system cost per unit area-time including facility set-up, outbound delivery and inbound replenishment cost per unit area-time. Obviously, the optimal number of facilities $N\rightarrow\infty$; otherwise, the average outbound cost goes to infinity. Hence, the optimization problem can be expressed as follows,
\begin{align}
z = &\min\limits_{{\cal{N}},\left\{x_i\right\},\left\{{\cal{A}}_i\right\}} \lim_{\begin{subarray}{l} N\rightarrow\infty\end{subarray}}
{\frac{f}{\sum\limits_{i=1}^{N} {A_{i}}}\sum\limits_{i=1}^{N} {\left(1+\kappa\int_{{\cal{A}}_{i}}{ \|x-x_i\| d{x}}\right)}+\frac{\kappa r}{N} \sum\limits_{i=1}^{N} {l_i}} \label{original_obj}\\
& \text {s.t. } \eqref{eq:A_i_def} -
\eqref{eq:l_i_def}.\notag\end{align}

\begin{figure}[!htbp]
 \centering
 \subfigure[Swath construction]{
 \label{Fig_Swath_Construction} 
 \includegraphics[width=0.3\textwidth]{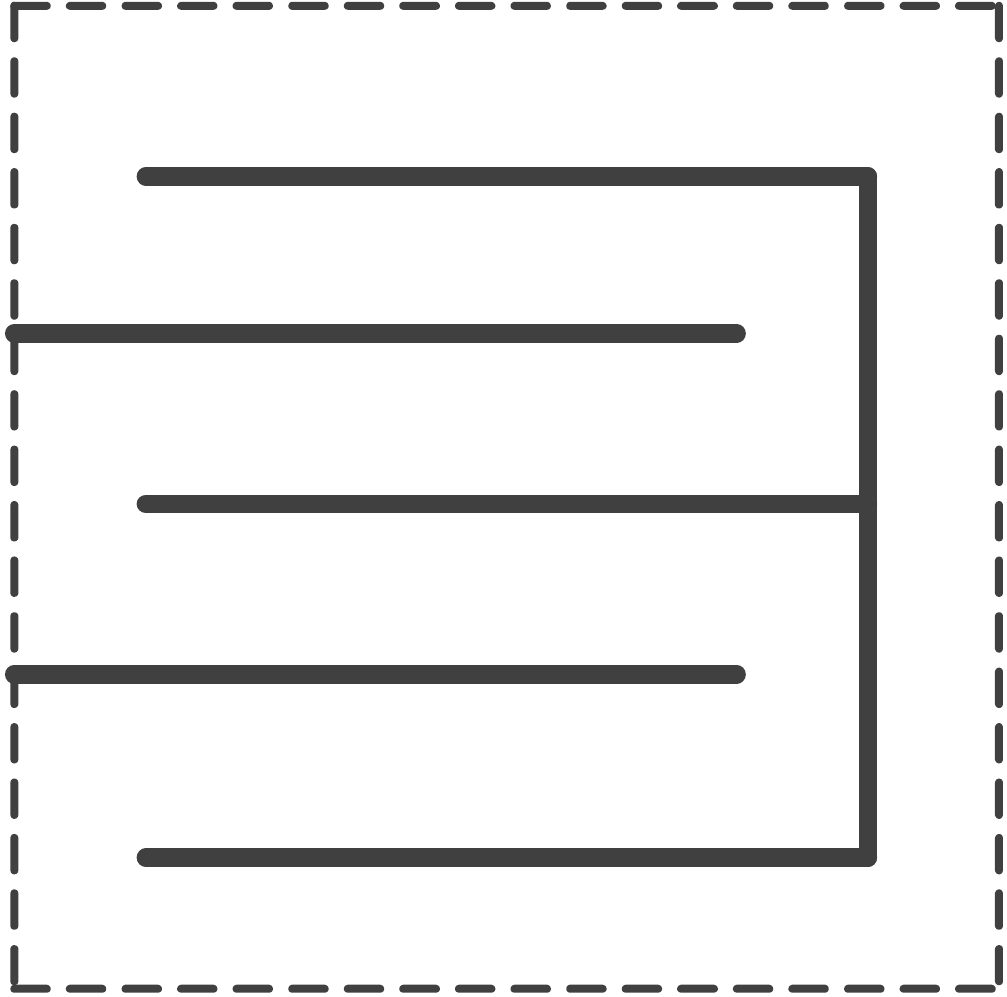}}
 \hspace{0.10in}
 \subfigure[Routing under Euclidean metric]{
 \label{Figure_Routing_Metric}
 \includegraphics[width=0.2\textwidth]{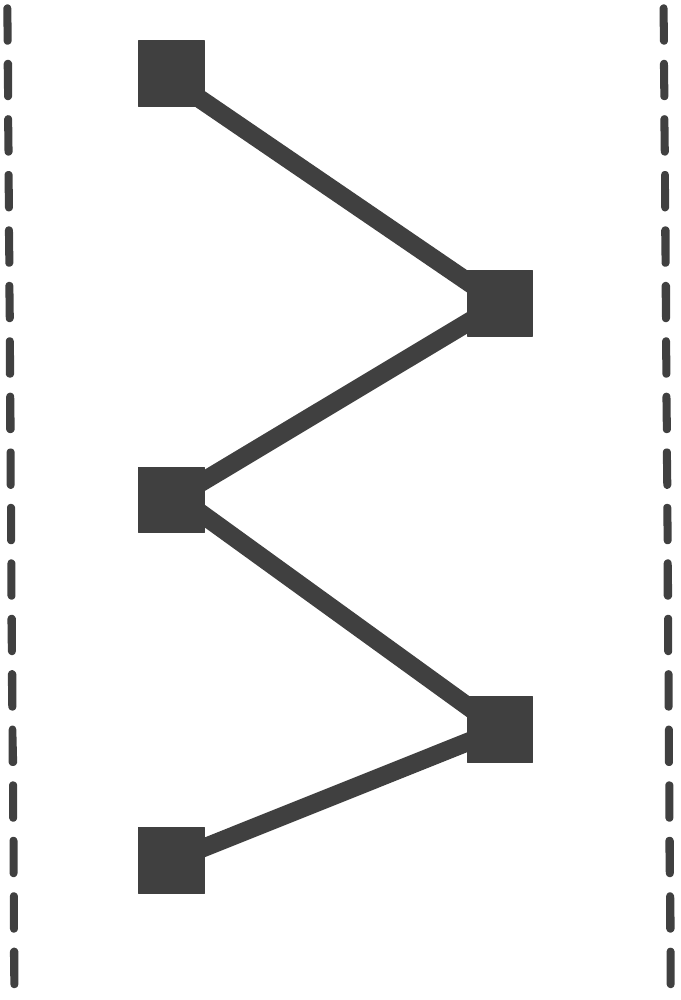}}
 \hspace{0.10in}
 \subfigure[Routing under $L_1$ metric]{
 \label{Figure_Routing_L_1}
 \includegraphics[width=0.2\textwidth]{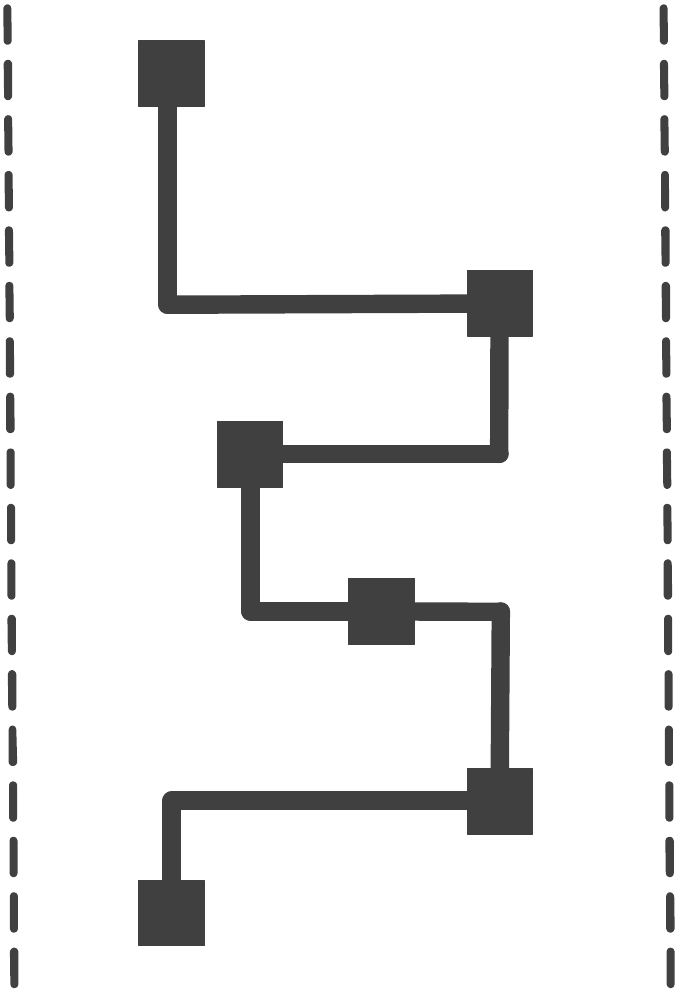}}
 \caption{Illustration of swath strategy for building asymptotically near-optimum TSP tour (Source: \cite{daganzo1984length})} \label{Figure_Daganzo_Tour}
\vspace*{-5pt}
\end{figure}

Only limited literature has addressed some simpler versions of this problem, mainly in the form of large-scale travelling salesman problem, planar facility location problem or uniform quantization problem. On the routing side, it is well-known \citep{beardwood1959shortest} that the optimal TSP tour length to visit $N$ randomly distributed points in an area of size $A$ asymptotically converges to $k\sqrt{NA}$, for some constant $k$, when $N \rightarrow \infty$. \citet{daganzo1984length} proposed a swath strategy for building asymptotically near-optimum TSP tour, as illustrated in Figure \ref{Figure_Daganzo_Tour}. The decision essentially reduces to determining the swath width to balance the trade-off between longitudinal travel (which is somewhat inversely proportional to the swath width), and local (or lateral) travel to reach customers (which increases monotonically with the swath width). These results hold for both Euclidean and rectilinear (i.e., $L_1$) metrics, and they shed light on the asymptotic behavior of our problem if we imagine the customers in each ${\cal A}_i, \forall i$ are spatially clustered at one point $x_i$ (instead of being spatially distributed). On the facility location side, \citet{Newell1973} pointed out that warehouse service regions should ideally have ``round'' shapes in order to minimize outbound delivery costs, although he also acknowledged the fact that round regions do not form a spatial partition. Similarly, \cite{Gersho79} conjectured that hexagonal shape is generally optimal in a two dimensional space. Later, \cite{Newman82} and \cite{haimovich1988extremum} respectively proved that under squared Euclidean metric and Euclidean metric, regular hexagonal service regions are optimal for outbound customer service. These results, albeit very relevant, do not solve our problem because they hold only when the inbound routing cost is ignored. In fact, \cite{geoffrion1979making} incorporated inbound freight cost into a warehousing location problem while approximating the service regions with identical squares, and also compared the effects of different cost components (fixed facility cost, outbound cost and inbound cost). Recently, \cite{Cachon11} compared three regular service region shapes (i.e., equilateral triangle, square and regular hexagon), showing that equilateral triangle tessellation is the best among these three when inbound transportation cost is relatively high. \cite{Carlsson13} further proposed several feasible tessellations, among which Archimedean spiral was proven to be asymptotically optimal under dominating inbound cost. This indicates that the consideration of inbound vehicle routing cost significantly affects the optimal spatial configuration of the transshipment system, making ``round'' service region shapes undesirable. Intuitively, inbound transportation cost tends to favor a design where facilities are clustered, which nevertheless twists the shape of the service region to become ``irregular''. However, to the best of our knowledge, no general results have been revealed regarding the optimal facility layout and service region configuration that achieve the best trade-off among the inbound, outbound, and facility costs.


Knowing the optimal shape of service regions would allow researchers to qualitatively approximate the outbound logistics cost in the facility service region, which significantly simplifies the modelling challenge. Hence, while the optimal shape was generally unknown, various studies had built their solution methods upon the ``likely'' optimal results. As early as the 1960s, \citet{Edwin64} compared four different market partition shapes (equilateral triangle, square, regular hexagon and circle) in a market equilibrium model to maximize each firm's profit within its market region; however, they did not prove that regular hexagon should be the optimal choice of market layout. While solving transshipment problems, \citet{Daganzo05} assumed round service regions to approximate the outbound cost at each transshipment level, which was later proven to provide a tight cost lower bound for the actual optimal design~\citep{Ouyang06}. \cite{Qi2010worst} developed several spatial partitioning approaches based on regular hexagonal heuristics to cover arbitrarily distributed demand points, and showed that the worst-case errors were bounded by a constant. \citet{Cui10} and \cite{Li10} integrated regular hexagonal shapes into the facility location design under probabilistic facility disruptions.


This paper aims to provide a rigorous foundation on the optimal or near-optimum spatial layout of transshipment facilities. We first present a new proof, different from that in ~\cite{haimovich1988extremum}, that the conclusion of \cite{Newman82} also holds for the Euclidean metric case, i.e., regular hexagon is still the optimal shape for facility service regions when inbound routing cost is negligible. Then, when inbound transportation cost is taken into consideration, we introduce a near-optimum shape -- a cyclic hexagon with two equal ``long'' sides and four equal ``short'' sides. We also provide formulas to compute the size and shape of the service region (including the length of TSP tour within each shape), the best facility layout, and the best total system costs for this near-optimum configuration. An infeasible lower bound is introduced by relaxation and idealization, which shows that the proposed cyclic hexagonal shape yields a very small gap. After all these discussion, we shift our focus to an $L_1$ metric plane and show that a cost lower bound can be achieved by a similar spatial configuration with elongated non-cyclic hexagons, which hence becomes exactly optimal, as long as they are properly oriented in the coordinate system. Numerical experiments are conducted to verify the correctness of our analytical results for both Euclidean and $L_1$ metrics. In so doing, we formulate a mixed-integer mathematical program to solve a discrete version of the transshipment location-routing problem. Finally, we conduct sensitivity analyses to reveal managerial insights, and discuss the potential impacts of inventory cost.

The remainder of the paper is organized as follows. Section~\ref{Without inbound truck} proves the optimal spatial configuration when inbound truck is ignored. Section~\ref{With inbound truck} derives cost upper bound and lower bound when inbound cost is non-negligible and compares several different spatial configurations. All the results in Sections ~\ref{Without inbound truck} and ~\ref{With inbound truck} hold for Euclidean metric. Section~\ref{sec_configuration_L_1} further discusses results for $L_1$ metric. Section~\ref{sec_Discussion} presents the numerical experiments, analyses, and discussion. Section~\ref{conclusion} concludes this paper.

\section{Gersho's Conjecture under Euclidean Metric: the Special Case with Negligible Inbound Cost}\label{Without inbound truck}

As a building block, we first show that regular hexagon service regions are optimal on a Euclidean plane if we consider only facility and outbound delivery costs. This result is commonly known as the Gersho's Conjecture \citeyearpar{Gersho79}. It was proven in \citet{haimovich1988extremum} by replacing each basic triangle with two right-angled triangles and then proved that the two right-angled triangles which share the same hypotenuse should be identical. After this, they derived the cost function for these two right-angled triangles and showed the convexity of this cost function, which eventually leads to regular hexagons. In this section, we provide a more concise proof. The basic idea is to derive a cost lower bound and then show such a lower bound is achieved by regular hexagons.

From Property~\ref{P1}, we can see that the boundary between any two adjacent facilities extends a triangle with either of the facility locations (e.g., $\Delta {EBD}$ and $\Delta {GBD}$ in Figure~\ref{Fig_Basic_Triangle}). Moreover, these two triangles obviously must be identical (i.e., $\Delta {EBD} \cong \Delta {GBD}$). For simplicity, we define the following:

\begin{definition}\label{Def_Basic_Triangle_Angle} Within each facility service region,
\begin{itemize}
\item A ``{\it{basic triangle}}'' is the one extended by the facility location and one side of the service region border; e.g., $\Delta {EBD}$ or $\Delta {GBD}$ in Figure~\ref{Fig_Basic_Triangle};
\item A ``{\it{basic angle}}'' is the internal angle of a basic triangle at the facility location; e.g., $\angle BED$ or $\angle BGD$ in Figure~\ref{Fig_Basic_Triangle}.
\end{itemize}
\end{definition}

In order to prove our main result, we first introduce the following lemma:
\begin{lemma}\label{L1_Given_area_angle} If a basic triangle has a fixed basic angle and a fixed area size, the isosceles shape minimizes the outbound cost.
\end{lemma}

\begin{figure}[htbp]
 \centering
 \subfigure[Original polygon]{
 \label{Fig_Original_Without} 
 \includegraphics[width=0.2\textwidth]{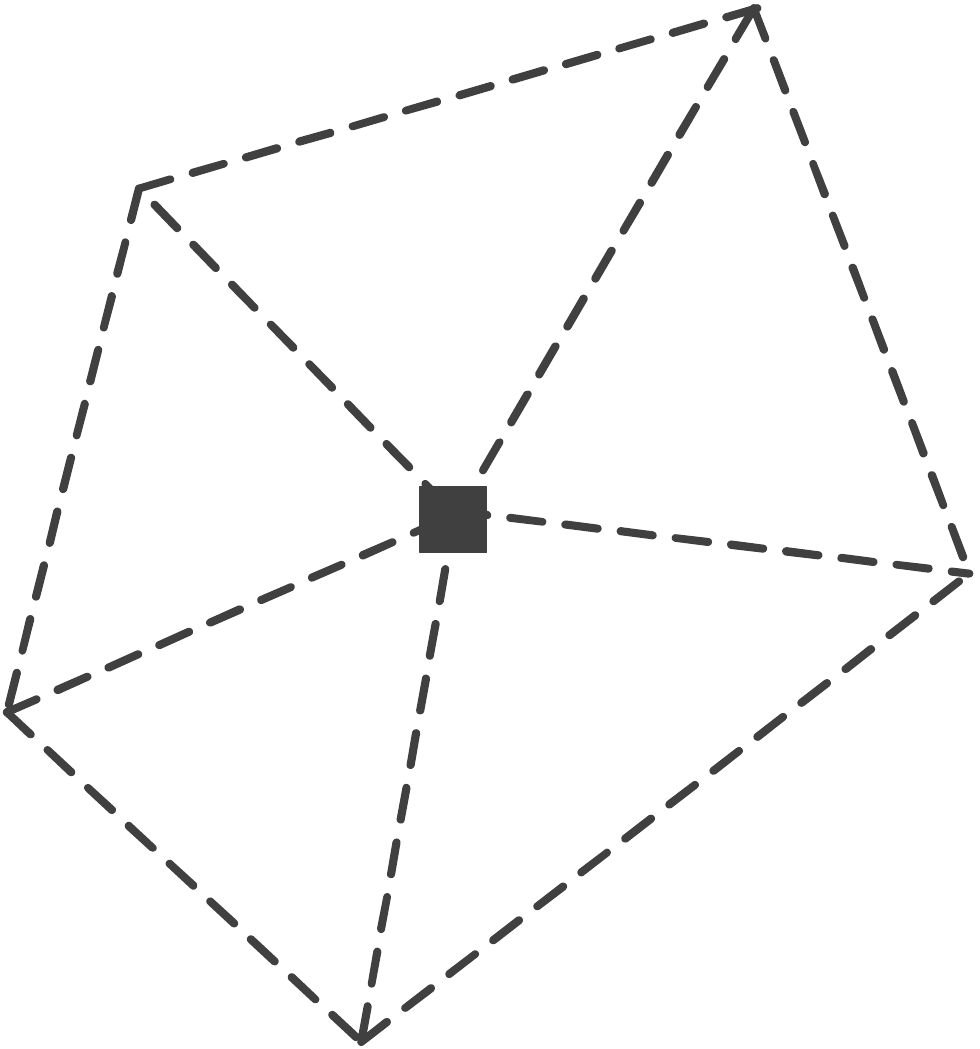}}
 \hspace{0.10in}
 \subfigure[Relaxed polygon]{
 \label{Fig_Relax_Without}
 \includegraphics[width=0.2\textwidth]{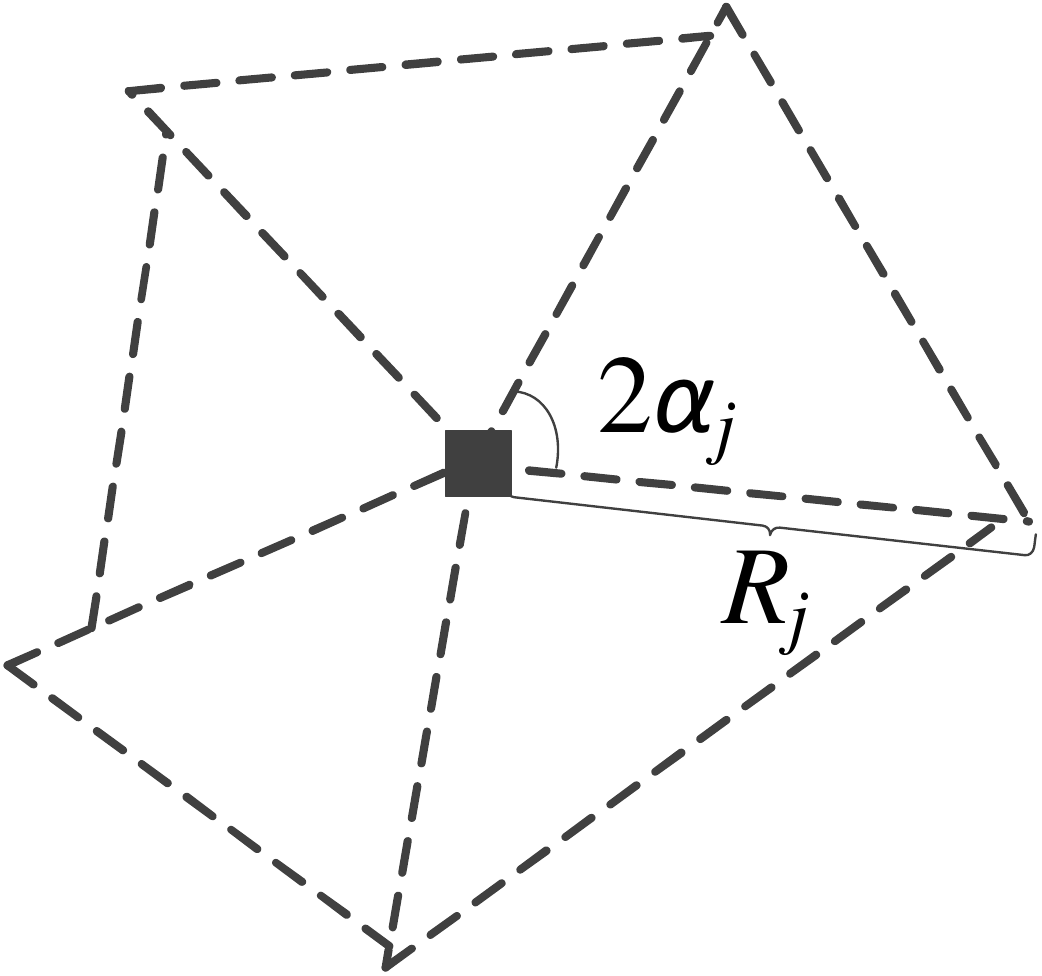}}
 \hspace{0.10in}
 \subfigure[Regular polygon]{
 \label{Fig_Regular_Without_n_5}
 \includegraphics[width=0.2\textwidth]{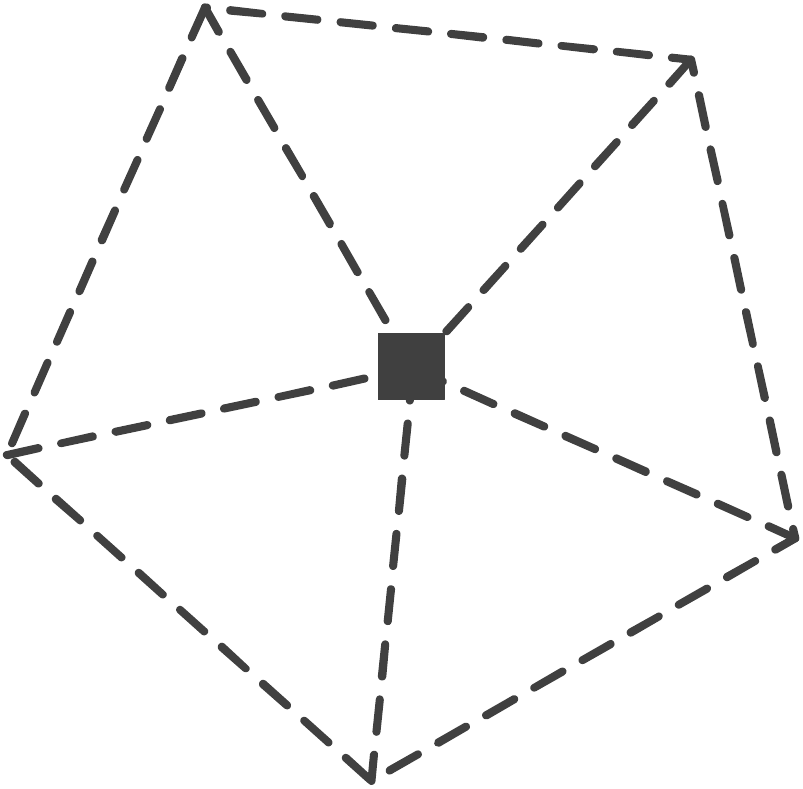}}
 \hspace{0.10in}
 \subfigure[Optimal polygon]{
 \label{Fig_Optimal_Without}
 \includegraphics[width=0.2\textwidth]{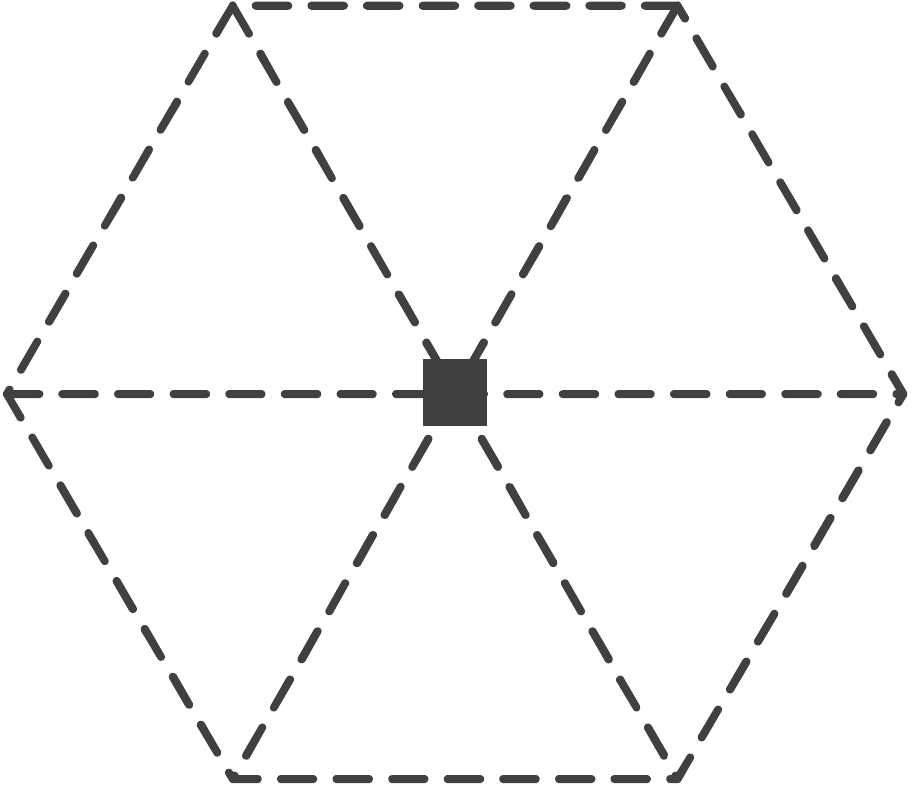}}
 \caption{Possible service region shapes}
 \label{verify_optimal} 
\end{figure}

Now consider a set of $n$ basic triangles (see Figure~\ref{Fig_Original_Without}) and the basic angles $2\alpha_1,\cdots,2\alpha_n$ satisfy $\sum_{j=1}^n\alpha_j=\theta,0<\theta\leq \pi$. If we relocate its customers so that each of these triangles becomes isosceles (with the same basic angle and area size), the resulting irregular area (see Figure~\ref{Fig_Relax_Without}) will have a lower total cost according to Lemma~\ref{L1_Given_area_angle}. Assume that the length of radial side of the $j$th isosceles basic triangle is $R_j$, then the outbound delivery cost for this service region is
\begin{equation}
 \sum\limits_{j=1}^n\frac{\kappa f R_j^3\cos^3\alpha_j}{3}\int_{-\alpha_j}^{\alpha_j}\frac{1}{\cos^3 t} dt\label{z_all}
\end{equation}
where $\sum\limits_{j=1}^nR_j^2\sin\alpha_j\cos\alpha_j=A$ and $\sum_{j=1}^n\alpha_j=\theta$.

Next, we can show that these triangles should be identical (see Figure~\ref{Fig_Regular_Without_n_5}) to further reduce cost. This is formally given in the following lemma.
\begin{lemma}\label{P_Given_Area_All}
When inbound cost is negligible, if total area size of $n$ basic triangles is $A$ and total angle degree is $2\theta$, all these basic triangles should be identical in order to minimize the outbound cost (i.e., $\alpha_j=\frac{\theta}{n},R_j=\left[\frac{A}{2n}\sin\frac{2\theta}{n}\right]^{-\frac{1}{2}},\forall j$).
\end{lemma}

Considering a facility service region with $\theta=\pi$, we have the following corollary.
\begin{corollary}\label{C_Given_Area_All}
When inbound cost is negligible, if the $i$th service region ${\cal A}_i$ has a fixed number of sides $n_i=n$ and area size $A_i=A$, the regular shape minimizes the outbound cost (i.e., $\alpha_j=\frac{\pi}{n},R_j=\left[\frac{n}{2A}\sin\frac{2\pi}{n}\right]^{-\frac{1}{2}},\forall j$) and the optimal facility and outbound cost is given by
\begin{equation}
f+ \frac{\kappa f A\sqrt{A}}{g(n)},\label{lb_z_n}
\end{equation}
where
\begin{equation}
g(n)=3n^{\frac{1}{2}}\left(\tan{\frac{\pi}{n}}\right)^{\frac{3}{2}}\left(\log\tan\left(\frac{\pi}{2n}+\frac{\pi}{4}\right)+\frac{\tan \frac{\pi}{n}}{\cos \frac{\pi}{n}}\right)^{-1}.\label{g_n_overall}
\end{equation}
\end{corollary}

%
%
%


Standard algebra shows that $\frac{d }{dx}[g^2(x)]>0$ and $\frac{d^2 }{dx^2}[g^2(x)]>0$, where $x\geq 3$. Hence, we have the following lemma.
\begin{lemma}\label{Lem_Function_G}
Function $g^2(x)$ is strictly concave and monotonically increasing with $x\in[3,\infty)$.
\end{lemma}


Now we consider the case of multiple facilities. Let $A$ denote the total area occupied by $N$ service regions; i.e.,
\begin{equation}
\sum\limits_{i=1}^{N} {A_i}=A. \label{Area_con_system}
\end{equation}

Since cost lower bound~\eqref{lb_z_n} holds for each facility $i\in {\cal N}$, the average facility and outbound cost per unit area-time satisfies
\begin{equation}
 z(N,A) \geq \frac{Nf}{A}+\sum\limits_{i=1}^{N} {\frac{\kappa f A_i\sqrt{A_i}}{Ag(n_i)}} \label{total_cost_system}
\end{equation}

The following lemma shows that the right hand side of~\eqref{total_cost_system}, while subject to \eqref{Area_con_system}, also has a lower bound.

\begin{lemma}\label{Lem_lB_Z_N_A}
$\sum_{i=1}^{N} {\frac{A_i\sqrt{A_i}}{Ag(n_i)}}\geq \frac{A^{\frac{1}{2}}}{\sqrt{N}g\left(\sum_{i=1}^N n_i/N\right)},\forall \{n_i\},N$, when $\sum_{i=1}^{N} {A_i}=A$, and equality holds only if $n_i=constant$ for all $i$.
\end{lemma}

Note that~\cite{Newman82} proved that $\sum_{i=1}^{N}n_i\leq{6N}$ since $\cal N$ is the a usual tessellation of the plane. Thus, we further have
\begin{equation}
 z(N,A) \geq \frac{Nf}{A}+\frac{\kappa f A^{\frac{1}{2}}}{\sqrt{N}g\left(\sum_{i=1}^N n_i/N\right)}
 \geq \frac{Nf}{A}+\frac{\kappa f A^{\frac{1}{2}}}{\sqrt{N}g\left(6\right)}\geq 3\sqrt[3]{\frac{\kappa^2f^3}{4g^2(6)}}.\label{total_cost_system_opt_n}
\end{equation}
The first inequality holds from Lemma~\ref{Lem_lB_Z_N_A}; the second inequality holds from Lemma~\ref{Lem_Function_G}; the third inequality holds by setting the optimal $A/N$ value to be $\left(\frac{\kappa}{2}\right)^{-\frac{2}{3}}\left(g(6)\right)^{\frac{2}{3}}$.

The final lower bound, i.e., the last term in \eqref{total_cost_system_opt_n}, turns out to be feasible, and hence optimal, as it can be achieved when $n_i=6,\forall i,$ and $A_i=\frac{A}{N}=\left(\frac{\kappa}{2}\right)^{-\frac{2}{3}}\left(g(6)\right)^{\frac{2}{3}}$. This implies that identical regular hexagon is the optimal shape for facility service regions (shown in Figure~\ref{Fig_Optimal_Without}).
\begin{proposition}\label{P_convex and decreasing}
When inbound cost is negligible, regular hexagon is the optimal shape of facility service region under Euclidean metric.
\end{proposition}
%

\section{Spatial Configuration under Euclidean Metric and Non-Negligible Inbound Cost}\label{With inbound truck}
In this section, we further consider inbound transportation cost in addition to facility cost and outbound delivery cost. We will try to first obtain a cost upper bound by constructing a reasonable feasible solution. Then we will derive a cost lower bound based on relaxation and idealization. After that, we show that the gap between these bounds is quite small, and hence the proposed feasible solution is near-optimum.

\subsection{Upper Bound}

\begin{figure}[htbp]
 \centering
 \subfigure[Cyclic polygons]{
 \label{Feasible_Solution}
 \includegraphics[width=0.4\textwidth]{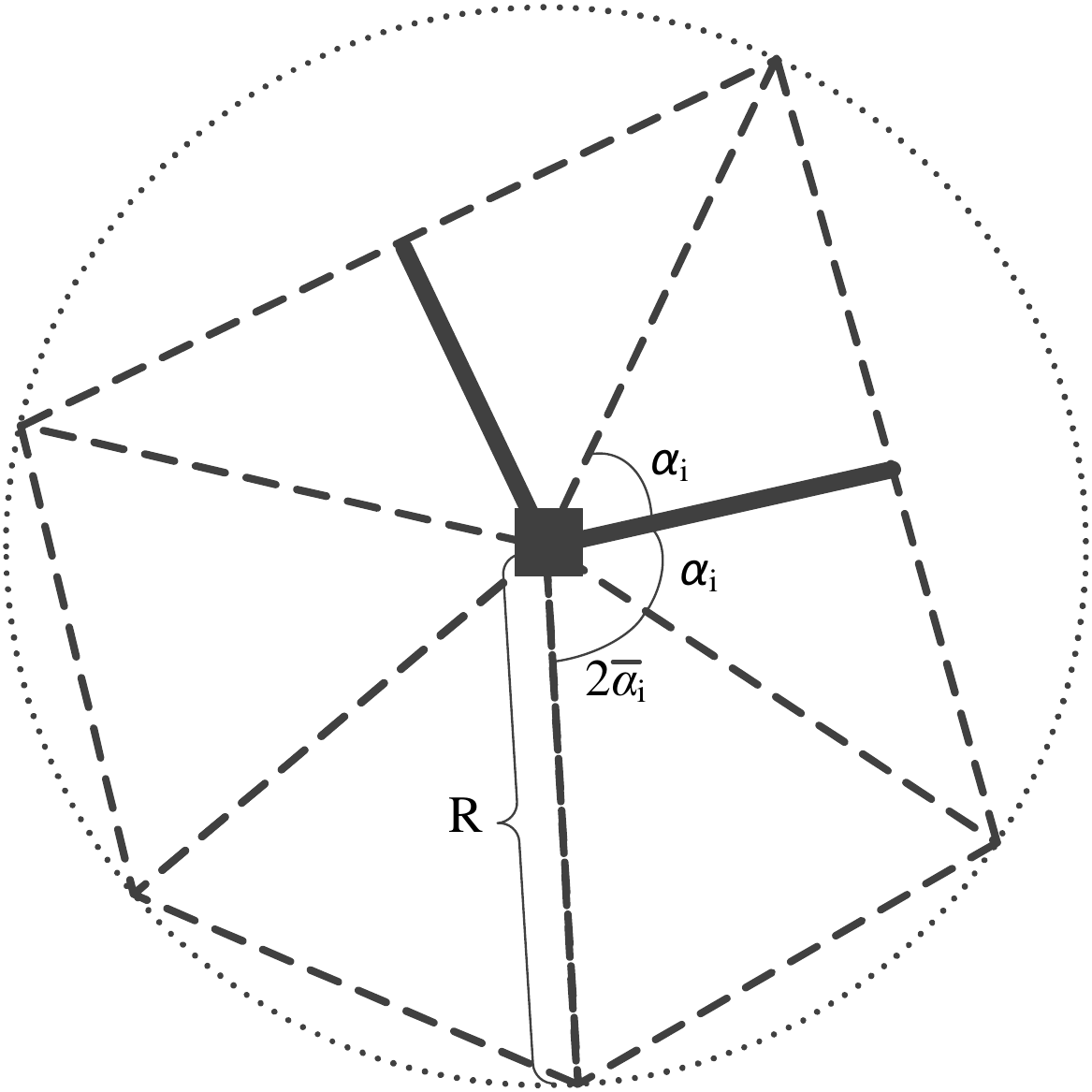}}
 \hspace{0.30in}
 \subfigure[A feasible tesselation]{
 \label{Fig_Optimal_With}
 \includegraphics[width=0.35\textwidth]{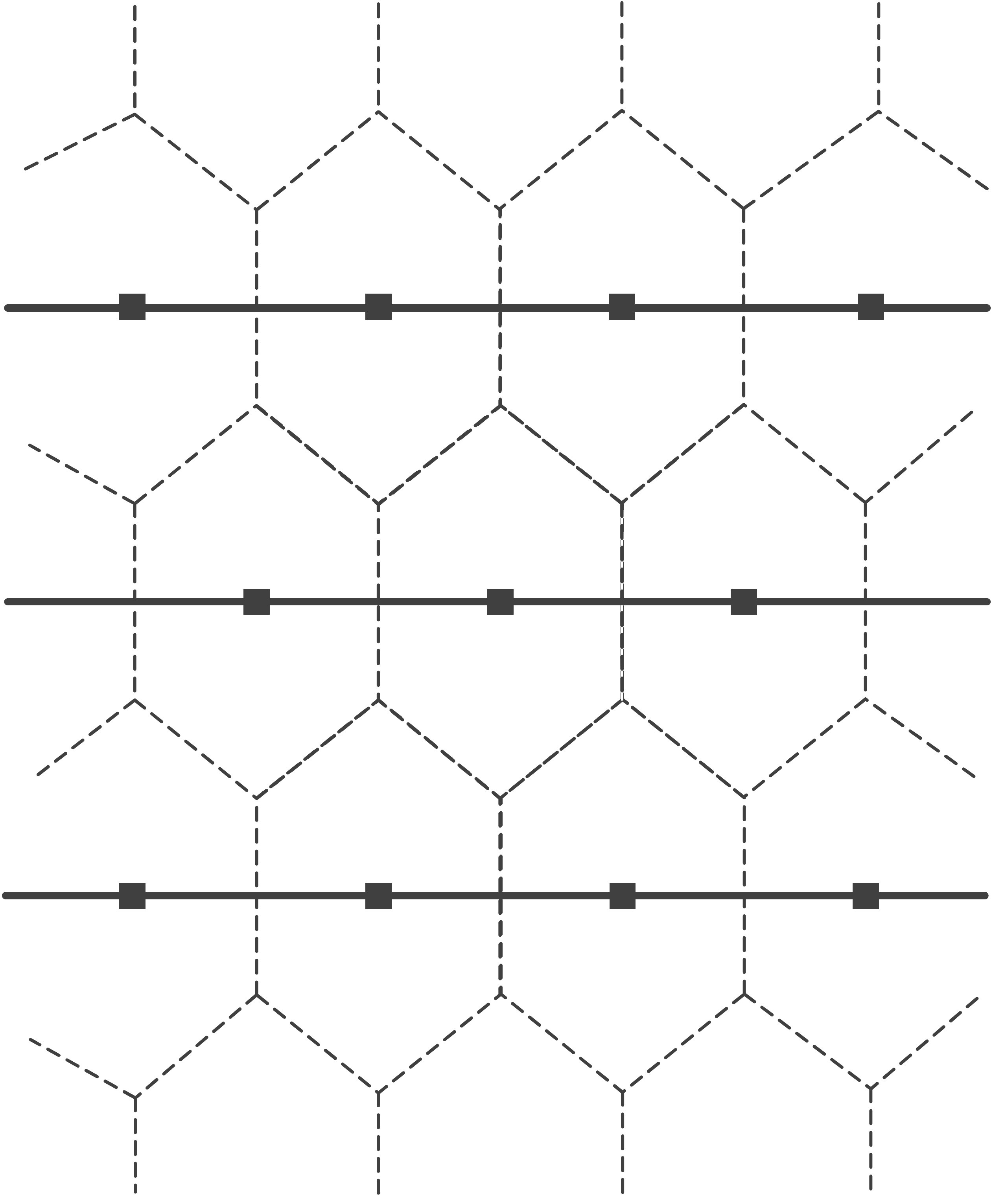}}
 \caption{A feasible tessellation and cost upper bound}
 \label{Fig_LB_UB} 
 \end{figure}

To construct an upper bound to \eqref{original_obj}, we first consider a set of $N \geq 1$ facilities ${\cal N} =\{1, 2, \cdots, N\}$. Facility $i\in {\cal N}$ serves a convex polygon service region ${\cal A}_i$ with $n_i \geq 3$ sides and area size $A_i$. Each polygon contains $n_i-2$ identical isosceles basic triangles, each with basic angle $2\bar{\alpha}_i$, with outbound delivery only and 2 identical isosceles basic triangles, each with basic angle $2\alpha_i$, which are passed by the inbound truck. The radial sides of all these triangles have an equal length of $R_i$, such that the $n_i$-sided polygon is cyclic; i.e., it is circumscribed by a circle of radius $R_i$. See Figure~\ref{Feasible_Solution} for an illustration.
These two types of basic angles $2\bar{\alpha}_i$, $2\alpha_i$ satisfy $(n_i-2) \bar{\alpha}_i+2\alpha_i=\pi, \bar{\alpha}_i>0$ and $\alpha_i>0$. Since the inbound truck must travel through the shortest distance within each service region, we must have $\alpha_i \geq \bar{\alpha}_i $ which implies that $\frac{\pi}{2}>\alpha_i\geq \frac {\pi}{n_i}$.

We shall be careful, that with this construct, the set of such cyclic polygons (even when $N \rightarrow \infty$) may or may not yet be feasible (i.e., forming a non-overlapping partition of $\Re^2$). However, as we shall see later, the cost-minimizer among this type of cyclic polygons happens to be feasible. To see this, we note that the average inbound and outbound costs for such an $n_i$-sided cyclic polygon service region are
\begin{equation}
 \frac{1}{A_i}\left[\frac{1}{3}(n_i-2)\kappa f R_i^3\cos^3\bar{\alpha}_i\int_{-\bar{\alpha}_i}^{\bar{\alpha}_i}\cos^{-3} t dt+\frac{2}{3}\kappa f R_i^3\cos^3\alpha_i\int_{-\alpha_i}^{\alpha_i}\cos^{-3} t dt+2\kappa r f A_iR_i\cos\alpha_i\right],\label{z_all_inbound_outbound}
\end{equation}
subject to $\frac{n_i-2}{2} R_i^2\sin2\bar{\alpha}_i + R_i^2\sin2\alpha_i=A_i$ and $(n_i-2)\bar{\alpha}_i+2\alpha_i=\pi$.

Before deriving the optimal basic angles $\bar{\alpha}_i$ and $\alpha_i$ that minimize \eqref{z_all_inbound_outbound}, we introduce a function in the following lemma.
\begin{lemma}\label{C_circumscribe_inbound_root}
For any given $n_i \geq 3$ and $r \in [0,\infty)$, the implicit equation $H(n_i,r,\alpha)=0$ has one and only one root with respect to $\alpha$ in the domain $\left[\frac{\pi}{n_i},\frac{\pi}{2}\right)$, which we denote by $\alpha^*(n_i,r ) $, where
\begin{align}
 H(n_i,r,\alpha)&= \sin\alpha\cos^2 {\frac {\pi-2\alpha}{n_i-2}} \log\tan \left( \frac{\pi}{4}+{\frac {\pi-2\alpha}{2(n_i-2)}} \right)-\cos^2\alpha\sin {\frac {\pi-2\alpha}{n_i-2}}\log \tan \left( \frac{\pi}{4}+\frac{\alpha}{2} \right)\notag \\
 &-r \sin 2\alpha\sin {\frac {\pi-2\alpha}{n_i-2}} -r (n_i-2)\cos {\frac {\pi-2\alpha}{n_i-2}} \sin^2 {\frac {\pi-2\alpha}{n_i-2}}.\label{equation_alpha_1_2_V2_H}
\end{align}
\end{lemma}

Now we are ready to show the optimal shape of the cyclic polygons defined above.
\begin{lemma}\label{P_circumscribe_inbound}
If an arbitrary service region takes the shape of an $n_i$-sided cyclic polygon defined above and has a fixed area size $A_i$, then the cost function \eqref{z_all_inbound_outbound} is minimized when the basic angles and the radial side length take the following values:
\[\alpha_i=\alpha^*(n_i,r ),\bar{\alpha}_i=\frac{\pi-2\alpha^*(n_i,r )}{n_i-2},R_i=\sqrt{A_i}
\left[(n_i-2)\sin\bar{\alpha}_i\cos\bar{\alpha}_i+\sin2\alpha_i\right]^{-\frac{1}{2}}.\]
Moreover, the total cost for this service region becomes
\begin{equation}
f+ \frac{\kappa f A_i\sqrt{A_i}}{g(n_i,r )},\label{lb_z_n_ith}
\end{equation}
where
\begin{equation}
g(n_i,r )=\frac{3\sin\bar{\alpha}_i \left(\sin 2\alpha_i +(n_i-2)\cos {\bar{\alpha}_i } \sin {\bar{\alpha}_i }\right)^{\frac{1}{2}}}{\sin\bar{\alpha}_i +\cos^2 {\bar{\alpha}_i } \log\tan \left( \frac{\pi}{4}+{\frac {\bar{\alpha}_i}{2}} \right)+4r \sin\bar{\alpha}_i\cos\alpha_i }.\label{g_n_inbound}
\end{equation}
\end{lemma}

Now we consider all $N$ such $n_i$-sided polygons, $i \in {\cal N} =\{1, 2, \cdots, N\}$, each with the optimal shape (e.g., basic angle $\alpha_i$ and radius $R_i$). Again, let $A = \sum_i^N A_i$ denotes the total area occupied by all $N$ service regions (i.e., \eqref{Area_con_system} holds). Since cost formula \eqref{lb_z_n_ith} holds for each service polygon $i \in {\cal N}$, the average cost per unit area across these $N$ service regions becomes:
\begin{equation}
 \frac{Nf}{A}+\sum\limits_{i=1}^{N} {\frac{\kappa f A_i^{\frac{3}{2}}}{Ag(n_i,r )}}. \label{total_cost_system_inbound}
\end{equation}

Proposition~\ref{Prop_lB_Z_N_A_inbound} below shows that, among all those shapes that form a spatial partition, \eqref{total_cost_system_inbound} reaches a minimum value when $n_i=6,\forall i$.
\begin{proposition}\label{Prop_lB_Z_N_A_inbound}
For polygons that form a partition of the Euclidean plane, for all $\{n_i\}, N \geq 1, r \geq 0$, and $A = \sum_{i=1}^{N} {A_i} $, we have $\sum_{i=1}^{N} {\frac{A_i^{\frac{3}{2}}}{Ag(n_i,r )}}\geq \frac{A^{\frac{1}{2}}}{\sqrt{N}g\left(\sum_{i=1}^N n_i/N,r \right)}\geq \frac{A^{\frac{1}{2}}}{\sqrt{N}g\left(6,r \right)} $. All equalities hold when $n_i=6$ for all $i$.
\end{proposition}

Hence, for all $r \geq 0$, the following holds.
\begin{equation}
 \eqref{total_cost_system_inbound} \geq \frac{Nf}{A}+\frac{\kappa f A^{\frac{1}{2}}}{\sqrt{N}g\left(\sum_{i=1}^N n_i/N,r \right)}
 \geq \frac{Nf}{A}+\frac{\kappa f A^{\frac{1}{2}}}{\sqrt{N}g\left(6,r \right)}\geq 3\sqrt[3]{\frac{\kappa^2f^3}{4g^2(6,r )}} = z_{ub}^*.\label{total_cost_system_opt_n_in}
\end{equation}
The last inequality becomes an equality only when $A_i=A/N, \forall i\in {\cal N}$.

In summary, the last term of \eqref{total_cost_system_opt_n_in} can be achieved when and only when $n_i=6,\forall i,$ and $A_i=\frac{A}{N}=\left(\frac{\kappa}{2}\right)^{-\frac{2}{3}}\left[g(6,r )\right]^{\frac{2}{3}}$. This implies that identical $6$-sided cyclic polygons (i.e., which we call ``{\it cyclic hexagons}'') is the cost minimizer among the class of cyclic polygons we have considered. Also, notice that $6$-sided cyclic polygons can obviously form a spatial partition, and hence the last term in \eqref{total_cost_system_opt_n_in} is achievable and feasible; see Figure~\ref{Fig_Optimal_With} for an illustration. Hence, $z_{ub}^*=3\sqrt[3]{\frac{\kappa^2f^3}{4g^2(6,r )}}$ is an upper bound, and likely a tight upper bound, of \eqref{original_obj}.

\subsection{Lower Bound}\label{sec_LB_L_2}

We now construct a cost lower bound which generalizes the asymptotic result in~\citet{Carlsson13}. Consider now the set of all solutions that incur a fixed inbound travel length $l$, a fixed number of facilities $N$ that collectively cover the customers in an area of total size $A$. The individual service regions in these solutions may take any shape, and some may not even be feasible if they do not form a spatial partition. The lowest possible total cost among these (relaxed) solutions will surely yield a cost lower bound.

\begin{figure}[ht]
 \centering
 \subfigure[]{
 \label{Figure_Isolated_Circle}
 \includegraphics[width=0.3\textwidth]{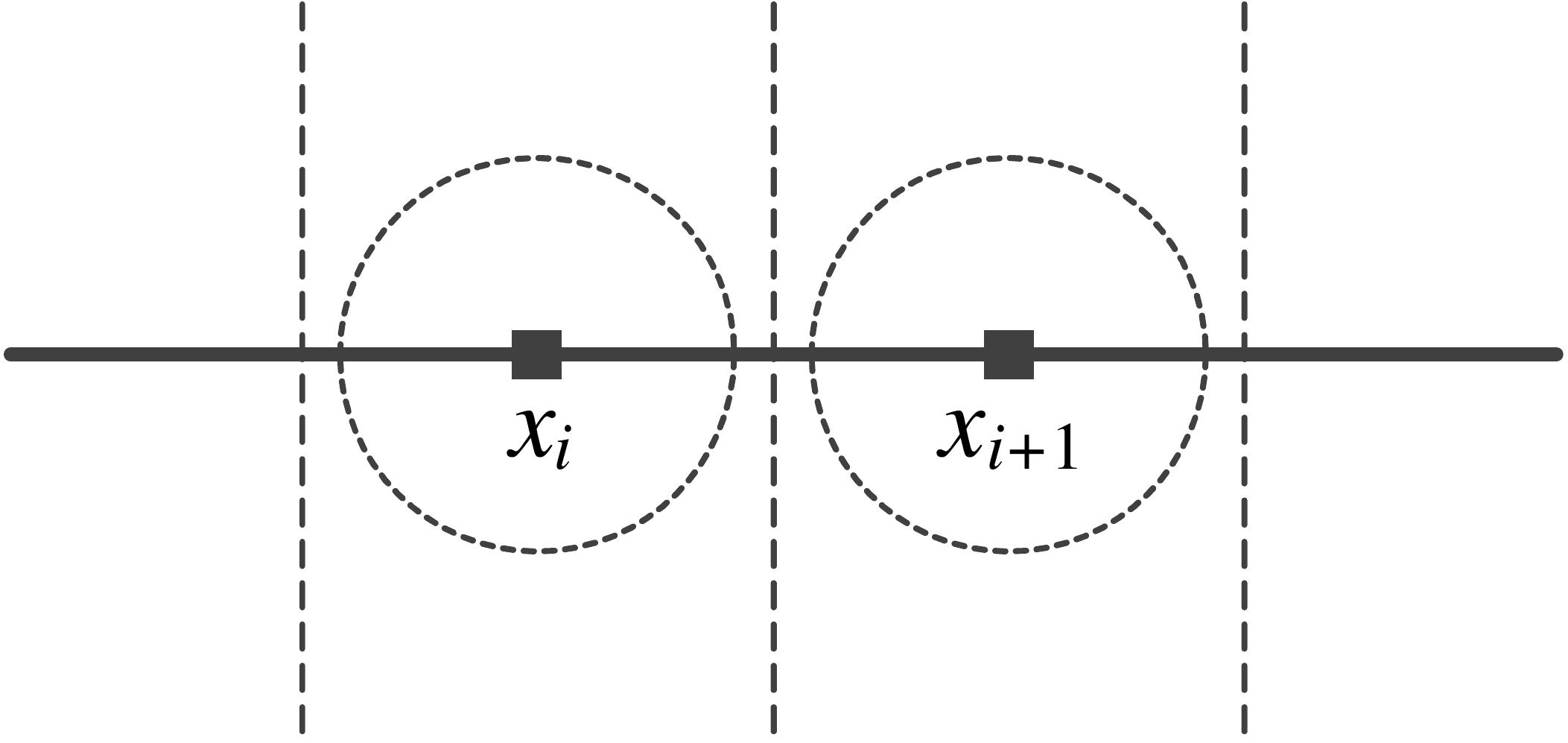}}
 \hspace{0.10in}
 \subfigure[]{
 \label{Figure_Boundary_Circle}
 \includegraphics[width=0.3\textwidth]{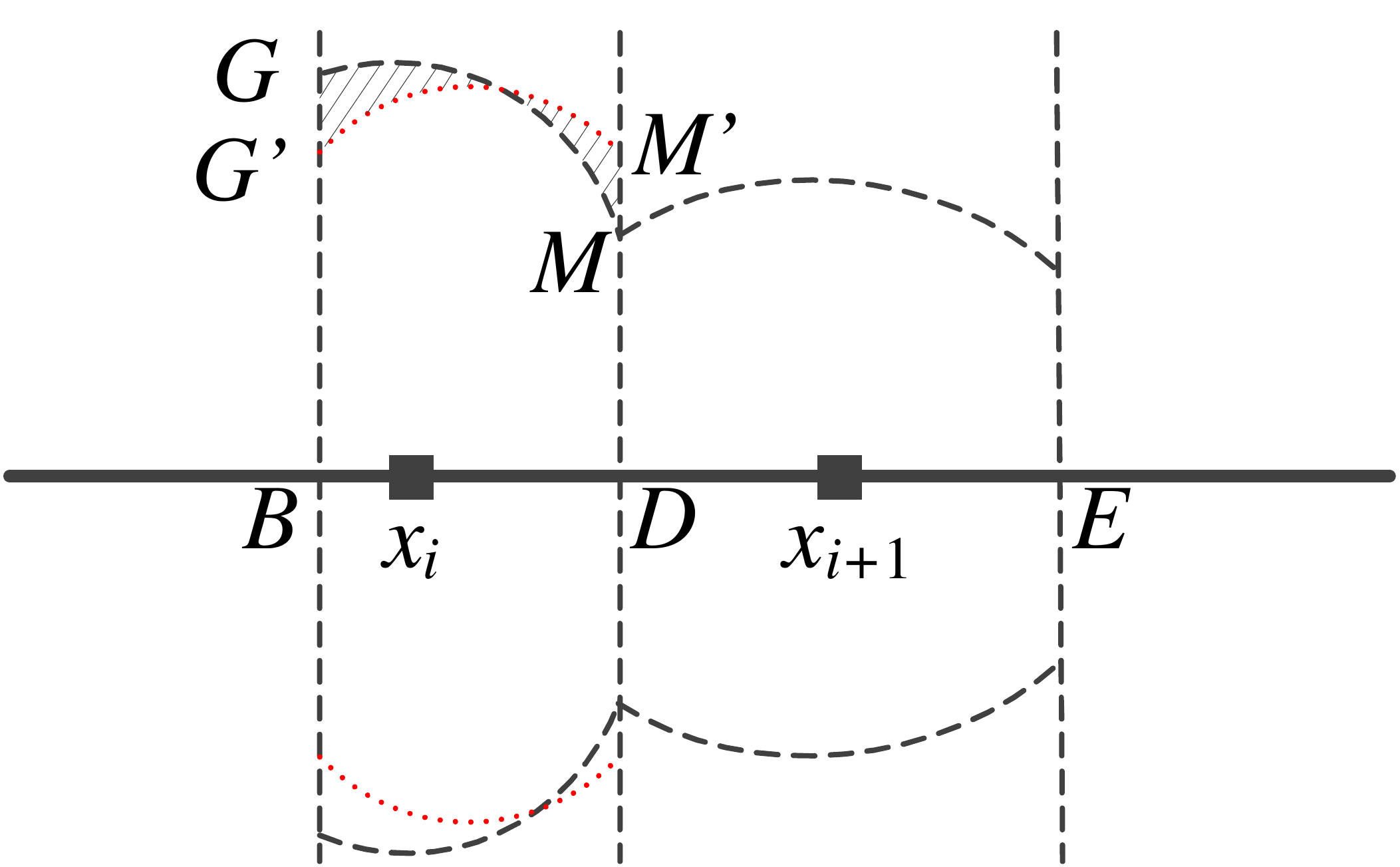}}
 \hspace{0.1in}
 \subfigure[]{
 \label{Figure_Line_Circle} 
 \includegraphics[width=0.3\textwidth]{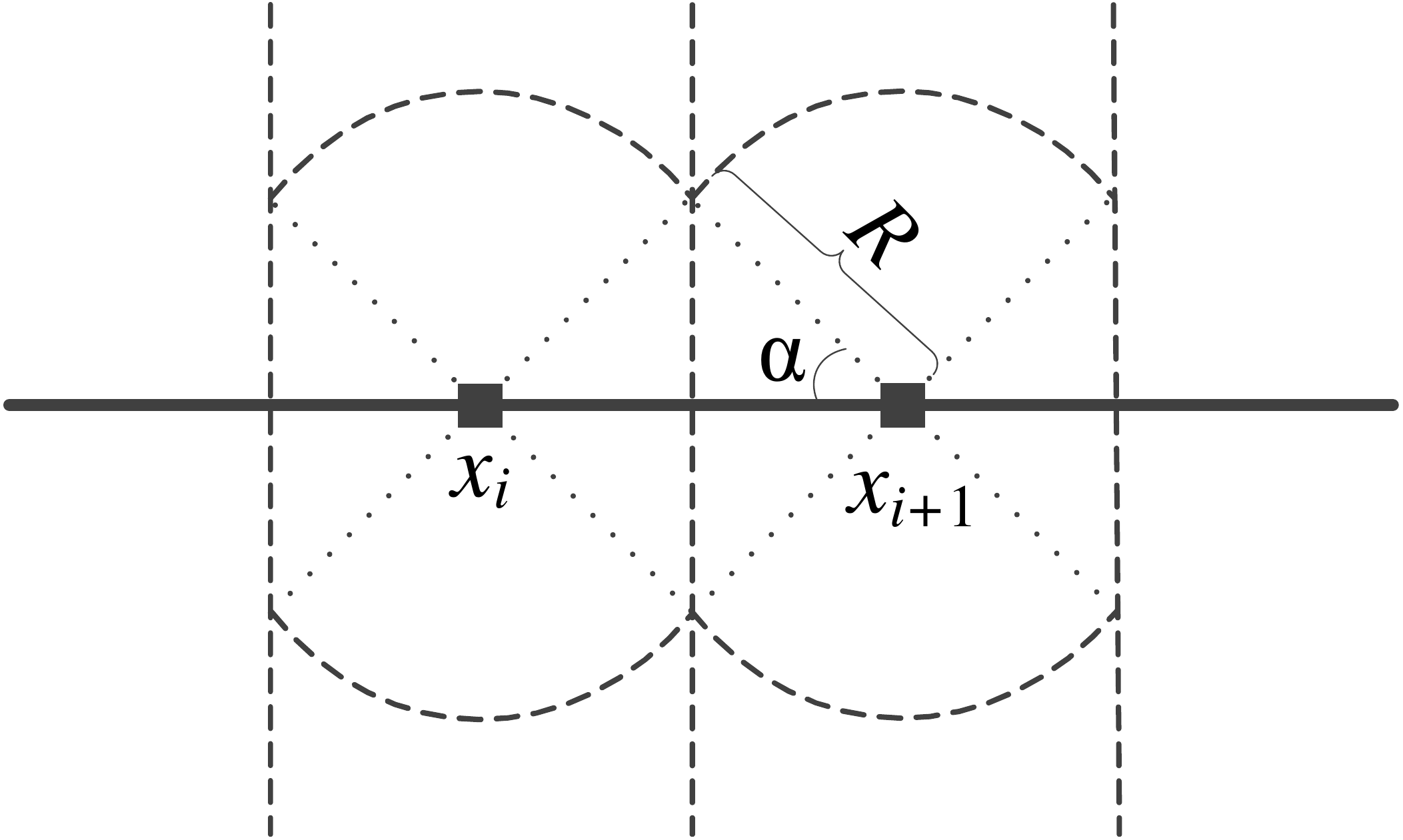}}
 \caption{Construction of a lower bound}
 \vspace*{-10pt}
 \label{Fig_Lem_Given_L_N_A_Boundary} 
\end{figure}

Note first that a circular shape minimizes the outbound cost~\citep{Ouyang06} for any given size of service regions. Thus, if the inbound travel length $l$ is larger than the total diameters of $N$ identical circles (each with area size $\frac{A}{N}$), then the case degrades to a trivial one where the optimal cost is achieved when all $N$ service regions take the shape of identical circles of radius $\left(\frac{A}{\pi N}\right)^{\frac{1}{2}}$. This case is illustrated in Figure~\ref{Figure_Isolated_Circle}. We shall note that this case never yields a good cost lower bound since the circular shape of the service regions will be far from forming a spatial partition, and that we can always shift the facility locations and their service regions along the TSP tour to reduce the length of inbound truck (without changing the total service region size $A$ nor the outbound costs).

Now we consider the more general case when $l\leq
2\left(\frac{NA}{\pi}\right)^{\frac{1}{2}}$. We first argue that the
lowest cost (for any fixed $l, N, A$) will be achieved when all $N$
facilities are along a straight line so as to minimize the potential
conflict among the facilities' customer sets. Otherwise, if there is an
``elbow'' facility along the TSP tour, we can always reduce the
outbound cost of this facility by straightening the corner, which
gives this facility the opportunity to serve more nearby customers,
while not changing the outbound customers of other facilities.

Recall from Property~\ref{P1} that the service region of each facility
should still be a Voronoi polygon. In case all facilities are along
a straight TSP line, two adjacent service regions
should be separated by a boundary line that is perpendicular to the TSP tour (see
Figure~\ref{Figure_Boundary_Circle}), and each facility
will only serve the customers within the two nearest boundaries. To minimize the
outbound cost of each facility within its boundaries, the optimal
service region should be the intersection of the area between the boundaries and
a circle, such that the maximum delivery distance is minimized. This
is easy to prove by contradiction, i.e., if this condition is not
satisfied, we can always trade customers from a farther location (e.g., those between $G$ and $G'$ in Figure~\ref{Figure_Boundary_Circle}) to nearer ones (e.g., those between $M$ and $M'$ in Figure~\ref{Figure_Boundary_Circle}) so as to minimize delivery cost while keeping the total service region size unchanged.

Finally, we will argue that the optimal locations and service regions of the facilities
should form a centroidal Voronoi tessellation, as shown in
Figure~\ref{Figure_Boundary_Circle}; i.e., any facility should also
be at the center of its service region. Otherwise, we can always
perturb the facility location within the fixed service region to
reduce the total outbound costs. Thus, according to Property~\ref{P1}, we must have $|x_iD|=|x_{i+1}D|$ in Figure~\ref{Figure_Boundary_Circle}; i.e., $|BD|=|DE|$. Hence, we can easily conclude that all
service regions must be identical, and all facilities are evenly
spaced along the straight TSP tour. This result is summarized in the
following lemma.
\begin{lemma}\label{Lem_Given_L_N_A_Boundary}
For a given $l$, $N$ and $A$ that satisfy $l\leq 2\left(\frac{NA}{\pi}\right)^{\frac{1}{2}}$, the lowest (per demand) system cost is achieved when (i) the TSP tour is a straight line along which all facilities lie evenly; (ii) the service regions of all facilities have the same size and shape; and (iii) each service region consists of two basic triangles and two pie shapes, as shown in Figure~\ref{Figure_Line_Circle}.
\end{lemma}
%

Intuitively, the service region in Figure~\ref{Figure_Line_Circle} is a special case of an $n_i$-sided cyclic polygon with $n_i\rightarrow\infty$. It shall yield a cost lower bound because it obviously cannot form a spatial partition. To obtain such a lower bound, we express $l$ as a function of $N,A$ and $\alpha$; i.e., $l=2\cos\alpha\left(AN\right)^{\frac{1}{2}}\left(\sin2\alpha +\pi-2\alpha\right)^{-\frac{1}{2}}$. Suppose the radius is $R$, then the lower bound can be achieved by solving the following minimization problem:
\begin{equation}
\min \frac{Nf}{A}+\frac{\kappa f N}{A}\left(\frac{2 }{3}R^3\cos^3\alpha\int_{-\alpha}^{\alpha}\cos^{-3} t dt+\frac{2}{3}\left(\pi-2\alpha\right)R^3+2r R\frac{A}{N}\cos\alpha\right),\label{eq_lower_bound}
\end{equation}
subject to $NR^2\left(\sin2\alpha +\pi-2\alpha\right)=A$. Note that for any given $A, N$, we can select the optimal $\alpha$ value (or the optimal shape of the service region).

The first order condition of \eqref{eq_lower_bound} with respect to $\alpha$ yields
\allowdisplaybreaks\begin{align*}
\frac{2\kappa f N\sin\alpha}{A\left(2\sin\alpha\cos\alpha+\pi-2\alpha\right)} \lim\limits_{n\rightarrow \infty}\left(n-2\right)H(n,r,\alpha)=0,
\end{align*}
where $H(n,r,\alpha)$ is defined in~\eqref{equation_alpha_1_2_V2_H}. Since $0< \alpha<\frac{\pi}{2}$, the above equation yields $\lim\limits_{n\rightarrow \infty}(n-2)H(n,r,\alpha)=0$, which has one and only one solution by Lemma~\ref{C_circumscribe_inbound_root}, $\alpha^*(\infty,r)$. Hence, when $\alpha=\alpha^*(\infty,r)$, we have
\allowdisplaybreaks\begin{align*}
\eqref{eq_lower_bound}\geq \frac{Nf}{A}+\frac{\kappa f A^{\frac{1}{2}}}{\sqrt{N}g\left(\infty,r \right)}
\geq 3\sqrt[3]{\frac{\kappa^2f^3}{4g^2(\infty,r )}},
\end{align*}
where $g(\infty,r ):=\lim\limits_{n\rightarrow \infty}g(n,r )$ and $g(n,r )$ is defined in~\eqref{g_n_inbound}. The last inequality becomes equality by choosing $\frac{A}{N}=\left(\frac{\kappa}{2}\right)^{-\frac{2}{3}}\left(g(\infty,r )\right)^{\frac{2}{3}}$.

The result is summarized in the following proposition.

\begin{proposition}\label{P_LB_inbound}
Let $\alpha^*(\infty,r)$ be the root of $\lim\limits_{n\rightarrow \infty}\left(n-2\right)H(n,r,\alpha) =0$ and $g(\infty,r ):=\lim\limits_{n\rightarrow \infty}g(n,r )$. A cost lower bound to \eqref{original_obj} under Euclidean metric is given by
\begin{equation}3\sqrt[3]{\frac{\kappa^2f^3}{4g^2(\infty,r )}}.\label{eq_LB_value}\end{equation}\end{proposition}

\subsection{Illustration: Impact of Service Region Shapes}\label{CATSF}

The upper and lower bounds presented in the previous subsections are general. To illustrate this, we compare the minimal costs of three intuitive shapes of service regions that can form a spatial partition: triangles ($n=3$), rectangles ($n=4$), and hexagons ($n=6$).

\begin{figure}[htbp]
 \centering
 \subfigure[Differences of optimal $\alpha$ between admissible shapes and the infinite polygon]{
 \label{Theta_difference} 
 \includegraphics[width=0.45\textwidth]{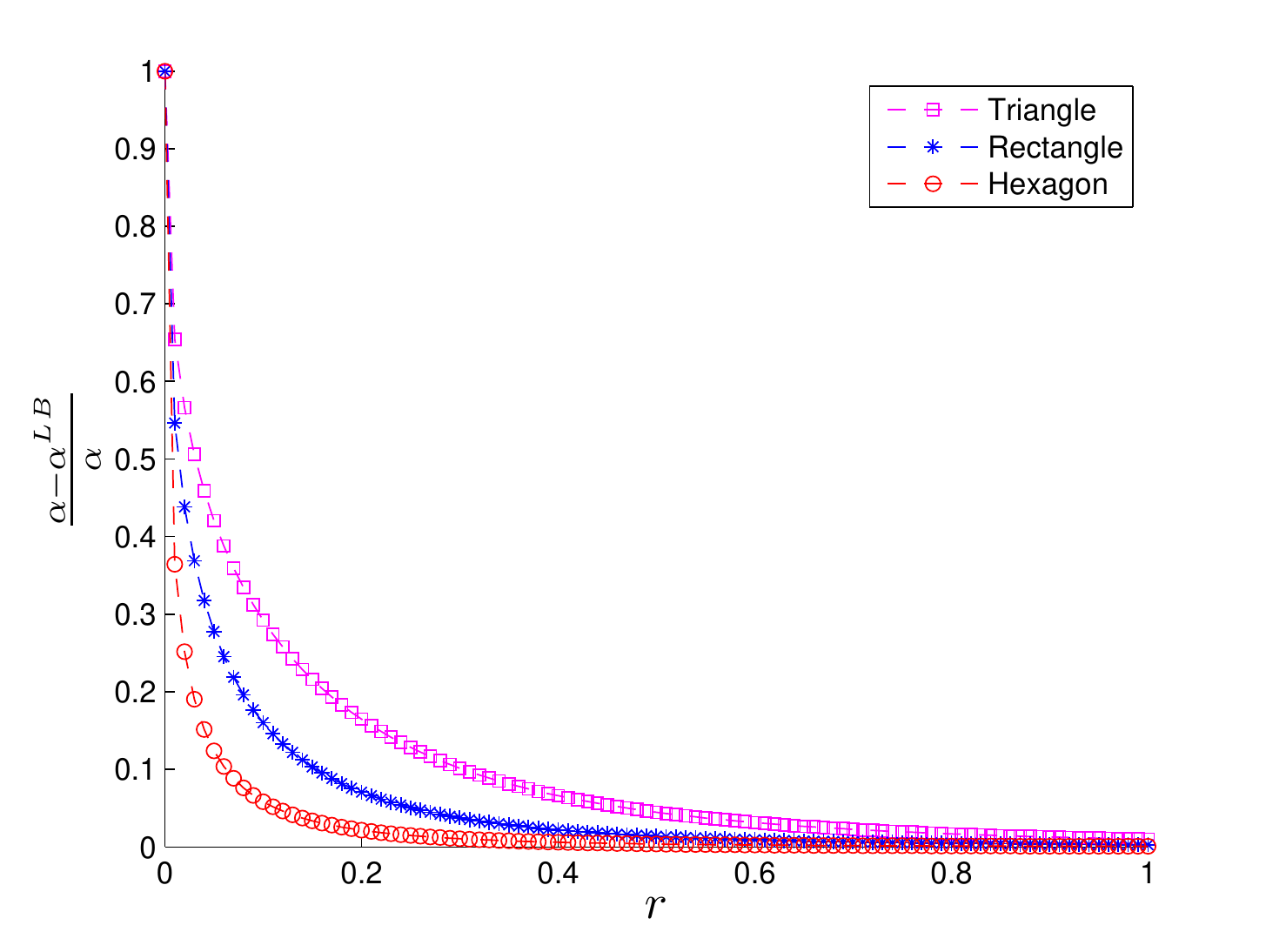}}
 \hspace{0.1in}
 \subfigure[Percentage cost differences between admissible shapes and the lower bound]{
 \label{Value_difference}
 \includegraphics[width=0.45\textwidth]{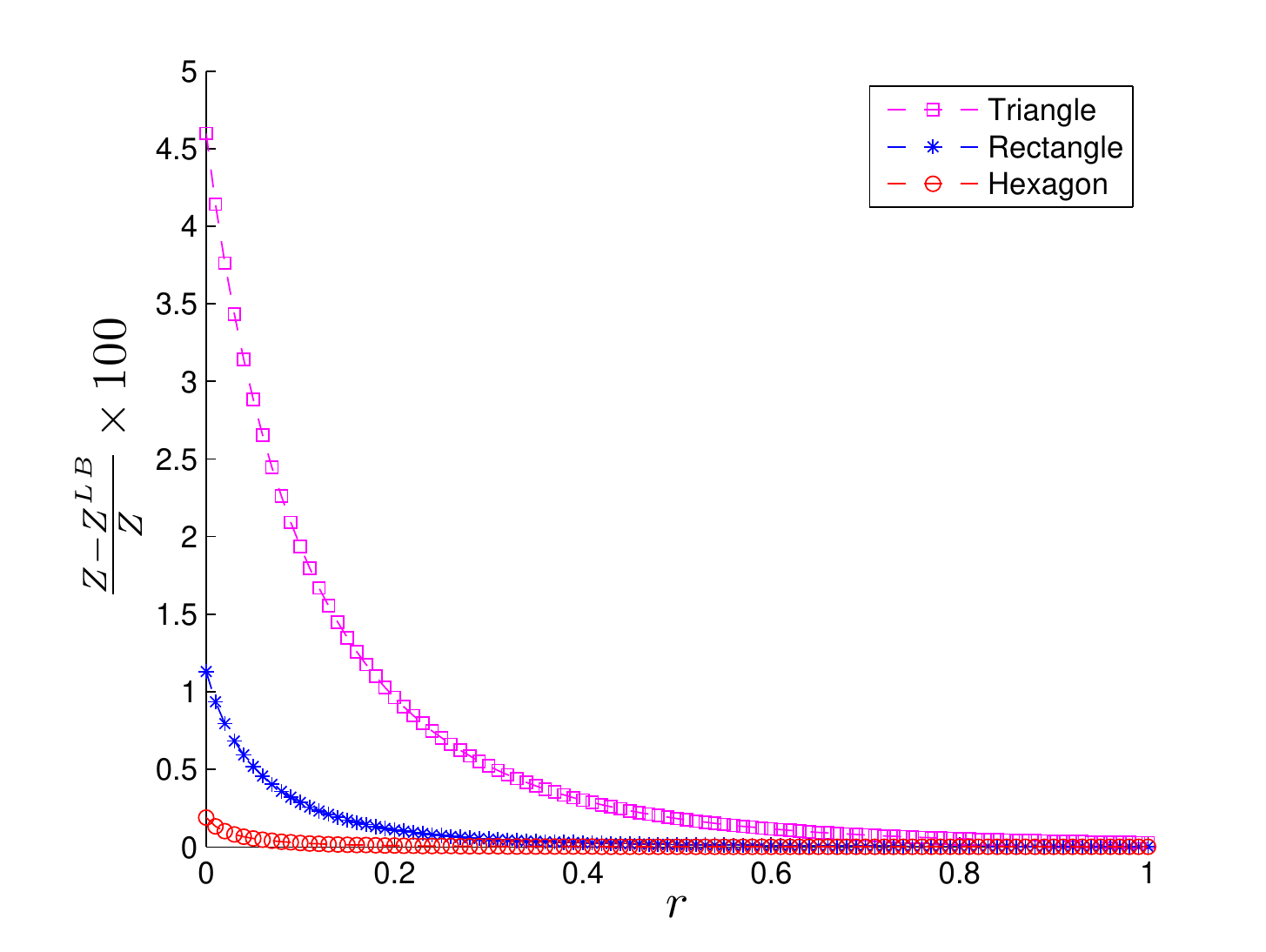}}
 \caption{Differences between admissible shapes and infinite polygon (an infeasible lower bound)}
 \vspace*{-5pt}
 \label{Differs} 
\end{figure}

For each shape, we fix $n_i=3,4,6,\forall i$, compute the optimal basic angles, and plug them into \eqref{total_cost_system_opt_n_in} (feasible cost upper bound) respectively.
The comparisons of the optimal $\alpha$ and optimal costs of these three special cases are shown in Figure~\ref{Differs} against the cost lower bound \eqref{eq_LB_value} for $n_i=\infty$. It can be seen that among these three intuitive shapes, the cyclic hexagon is the best.

It shall be noted that the differences of optimal angle $\alpha$ or optimal costs among these three shapes reduce as $r$ increases, so when $r$ is large enough, the spatial configuration of any one of these three shapes make no obvious difference. When $r \rightarrow 0$, the relative gap between the $n \rightarrow \infty$ lower bound and these feasible shapes grow. Fortunately, for the cyclic hexagon, the percentage cost gap remains quite small ($0.3\%$), suggesting strongly that the cyclic hexagons are near-optimum shapes.

\section{Spatial Configuration under $L_1$ Metric}\label{sec_configuration_L_1}

\subsection{Negligible Inbound Cost}\label{L_1_C_0}
\begin{figure}[ht]
 \centering
 \subfigure[Orignal service region]{
 \label{Figure_Original_L_1} 
 \includegraphics[width=0.3\textwidth]{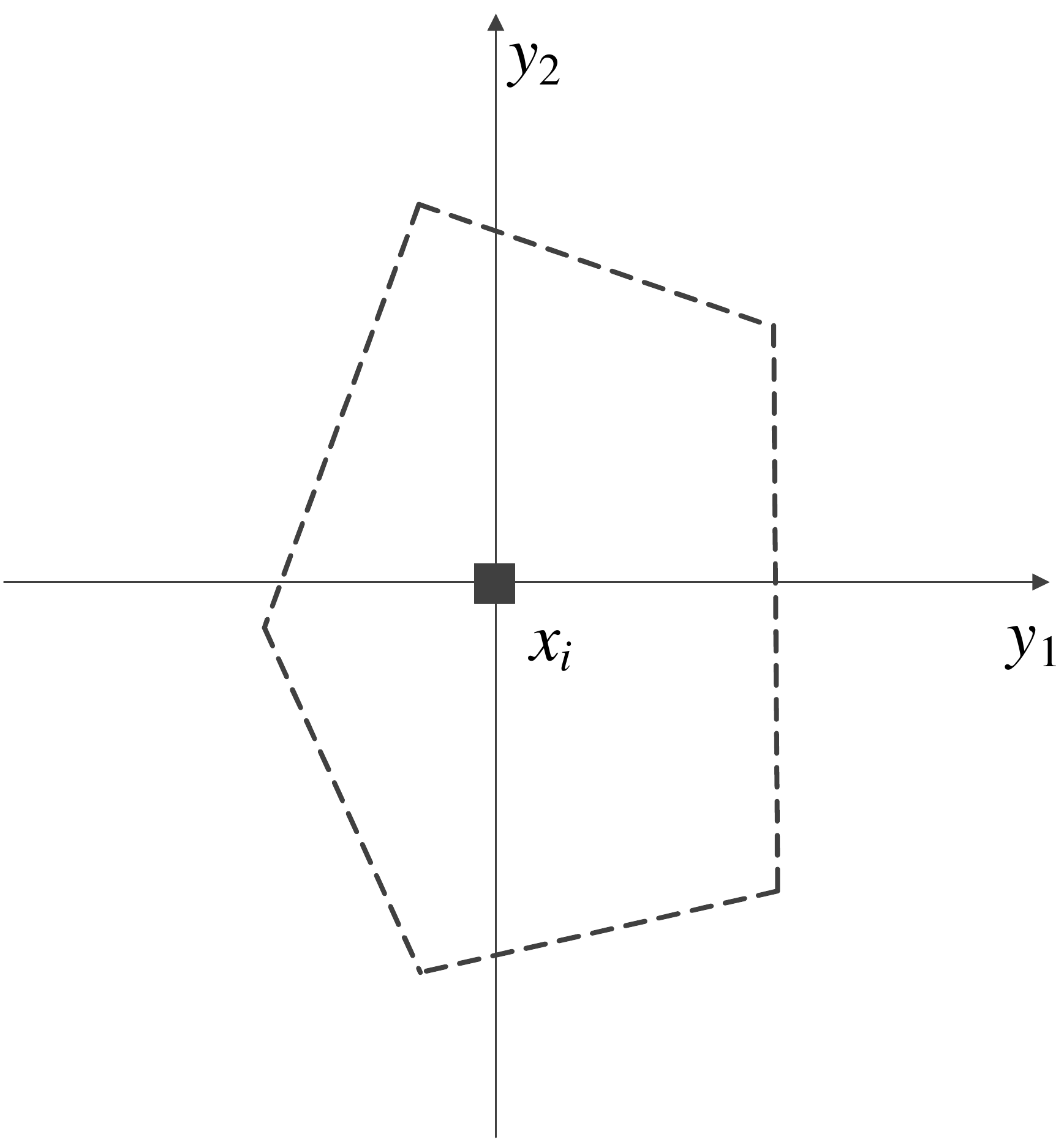}}
 \hspace{0.15in}
 \subfigure[Reshape the service region]{
 \label{Figure_Reshape_L_1}
 \includegraphics[width=0.3\textwidth]{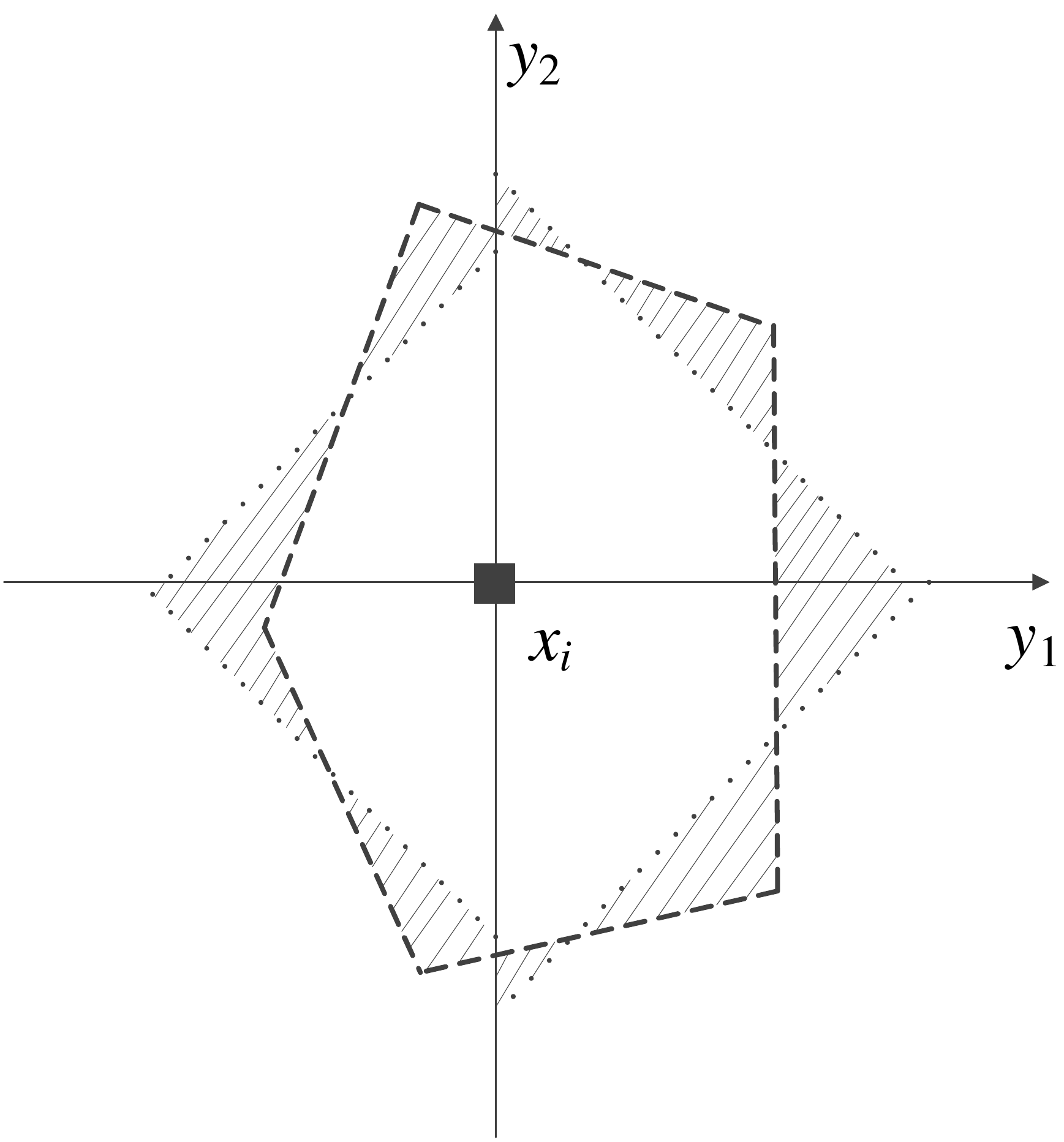}}
 \hspace{0.15in}
 \subfigure[Optimal shape]{
 \label{Figure_Square_L_1}
 \includegraphics[width=0.3\textwidth]{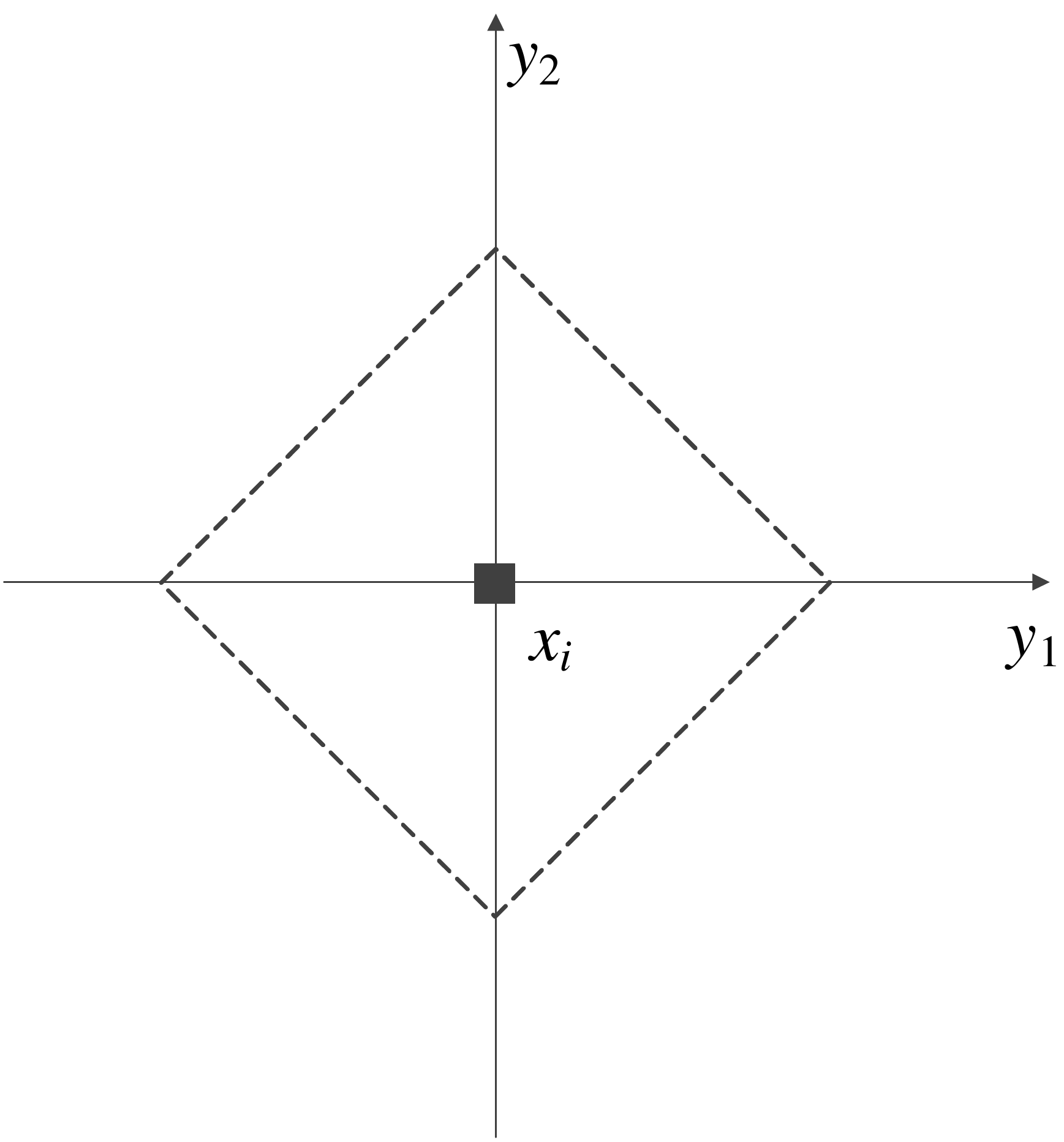}}
 \caption{Illustration of relaxation of a service region with $L_1$ metric}
 \vspace*{-5pt}
 \label{Fig_Relaxation_L_1} 
\end{figure}

Consider an arbitrary service region ${\cal A}_i$ on an $L_1$ metric plane. If we set the origin of the coordinate axes $y_1, y_2$ at the facility location $x_i$, the service region will be divided into four non-overlapping quadrant parts (i.e., ${\cal A}_{i1},{\cal A}_{i2},{\cal A}_{i3},{\cal A}_{i4}$), as shown in Figure~\ref{Figure_Original_L_1}. In each quadrant, there exists a line in the form of $|y_1|+|y_2|=$constant, such that the area of the resulting isosceles right-angled triangle equals that of the original quadrant part (see Figure~\ref{Figure_Reshape_L_1}). Note that all points on such a line have an equal travel distance to the facility. We shall easily see that the outbound service cost for the four isosceles right triangles is lower than that for the original ${\cal A}_i$, since our construct of these triangles can be done by simply re-locating some of the original customers to a nearer location (e.g., among the shaded areas in Figure~\ref{Figure_Reshape_L_1}). Following similar arguments of Lemma 2, we can also see that the outbound cost is further minimized when all four isosceles right triangles have the same size; i.e., ${ A}_{i1}={ A}_{i2}={A}_{i3}={ A}_{i4}=\frac{A_i}{4}$. As such, the triangles form a square, as shown in Figure~\ref{Figure_Square_L_1}), whose total outbound cost becomes
\begin{equation}
\frac{\sqrt{2}}{3}\kappa f A_i^{\frac{3}{2}}. \label{z_i_lowerbound_L_1}
\end{equation}

Now we consider a set of service regions ${\cal N}$ that form a partition, with a total area size equal to $A=\sum_{i \in {\cal N}} A_i$. Since cost lower bound~\eqref{z_i_lowerbound_L_1} holds for each region ${\cal A}_i, i\in {\cal N}$, the average facility and outbound cost per unit area-time satisfies the first equality below:
\begin{equation}
 z(N,A) \geq \frac{Nf}{A}+\frac{\sqrt{2}}{3}\kappa f \sum\limits_{i=1}^{N} {\frac{A_i^{\frac{3}{2}}}{A}}\geq\frac{Nf}{A}+\frac{\sqrt{2}}{3}\kappa f\left(\frac{A}{N}\right)^{\frac{1}{2}}
 \geq 3\sqrt[3]{\frac{\kappa^2f^3}{18}}. \label{total_cost_system_L_1}
\end{equation}
The second inequality obviously holds due to function convexity, and it becomes equality when $A_i=A/N, \forall i\in{\cal N}$; the third inequality becomes equality by setting the value of $A/N$ to $\left(\frac{\sqrt{2}\kappa}{6}\right)^{-\frac{2}{3}}$.

The final lower bound in \eqref{total_cost_system_L_1} is feasible and can be achieved when $n_i=4,$ and $A_i=\frac{A}{N}=\left(\frac{\sqrt{2}\kappa}{6}\right)^{-\frac{2}{3}}, \forall i \in {\cal N}$. This implies that identical square is the optimal shape for facility service regions under $L_1$ metric.
\begin{proposition}\label{P_Optimal_L_1}
When inbound cost is negligible, the optimal shape of facility service region under $L_1$ metric is square with diagonals parallel to the coordinate axes.
\end{proposition}

\subsection{Non-Negligible Inbound Cost}

\begin{figure}[ht]
 \centering
 \subfigure[]{
 \label{Figure_Zigzag_L_1}
 \includegraphics[width=0.3\textwidth]{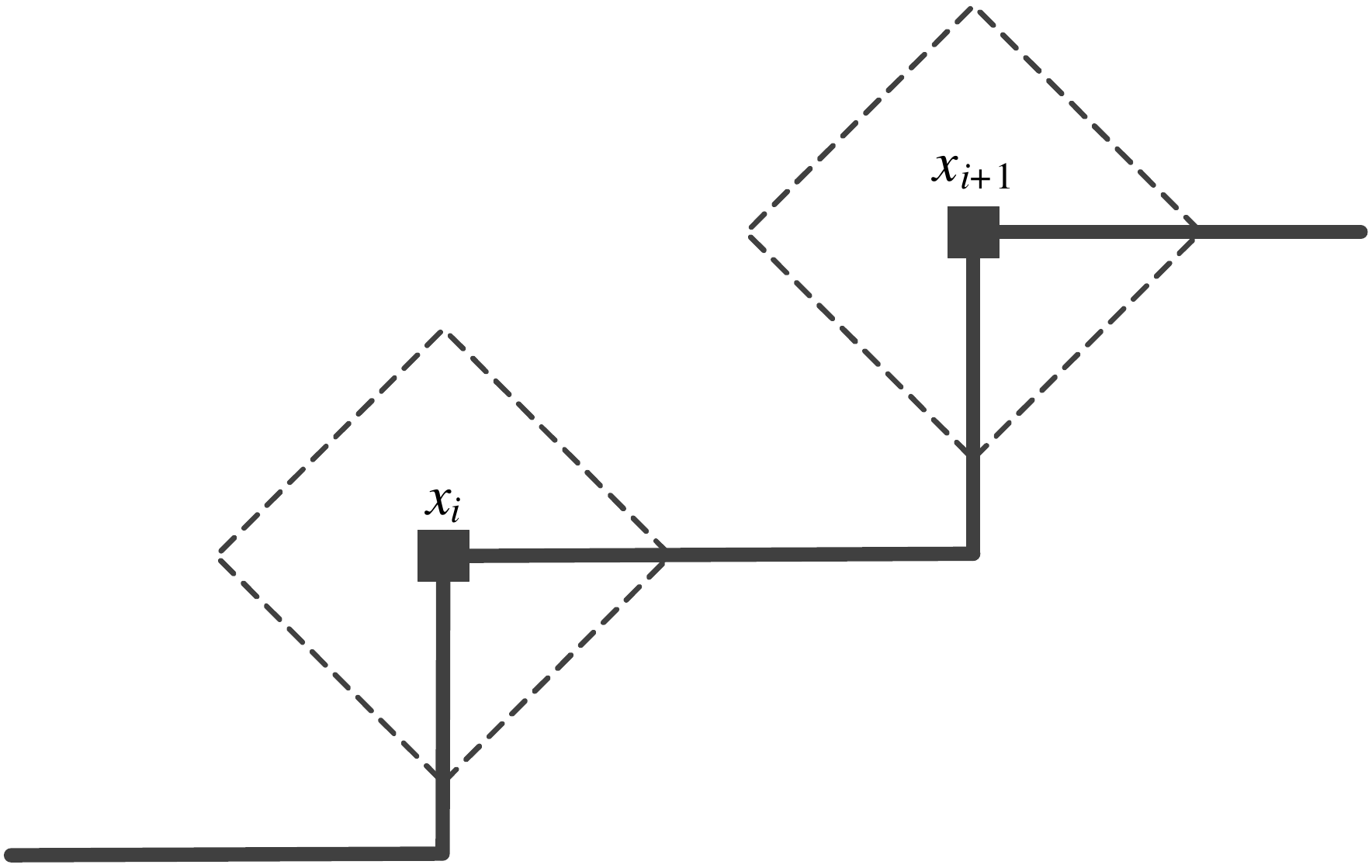}}
 \hspace{0.10in}
 \subfigure[]{
 \label{Figure_Zigzag_Overlap_L_1}
 \includegraphics[width=0.30\textwidth]{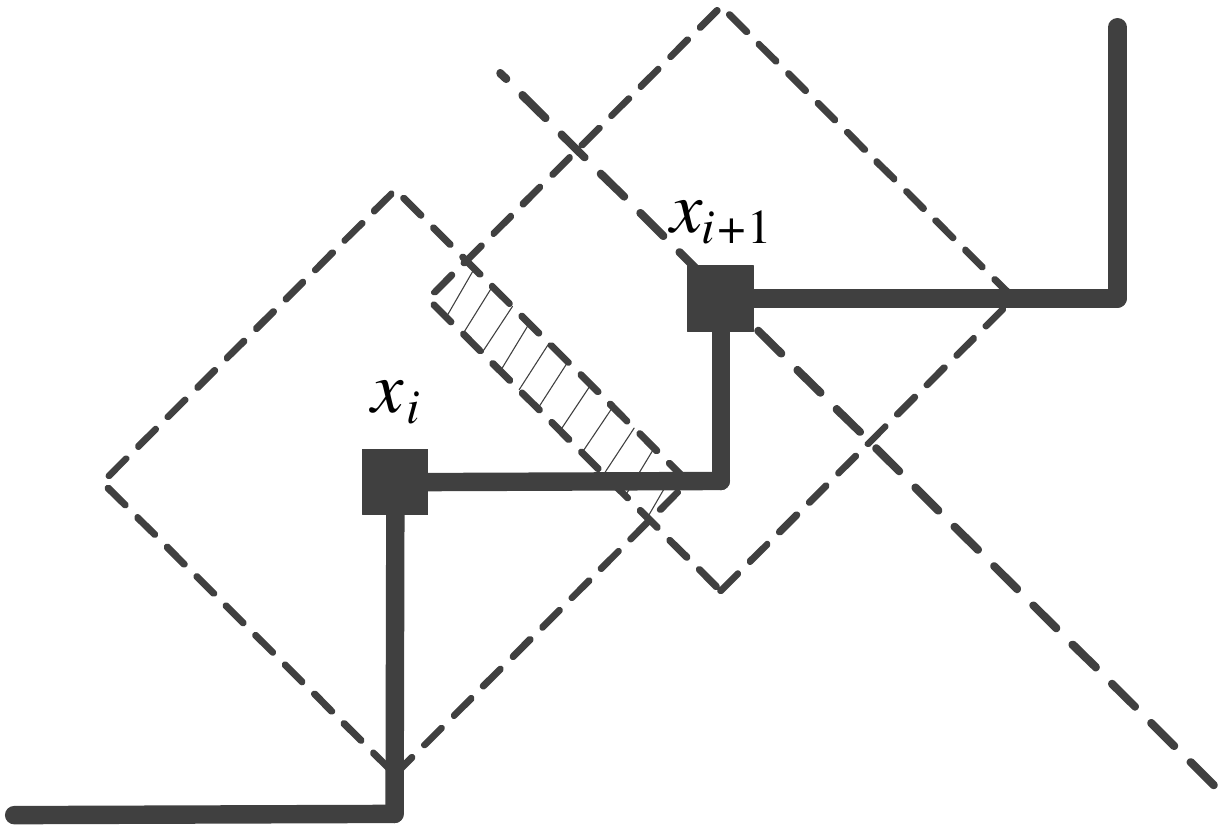}}
 \hspace{0.10in}
 \subfigure[]{
 \label{Figure_Zigzag_Overlap_Line_L_1}
 \includegraphics[width=0.3\textwidth]{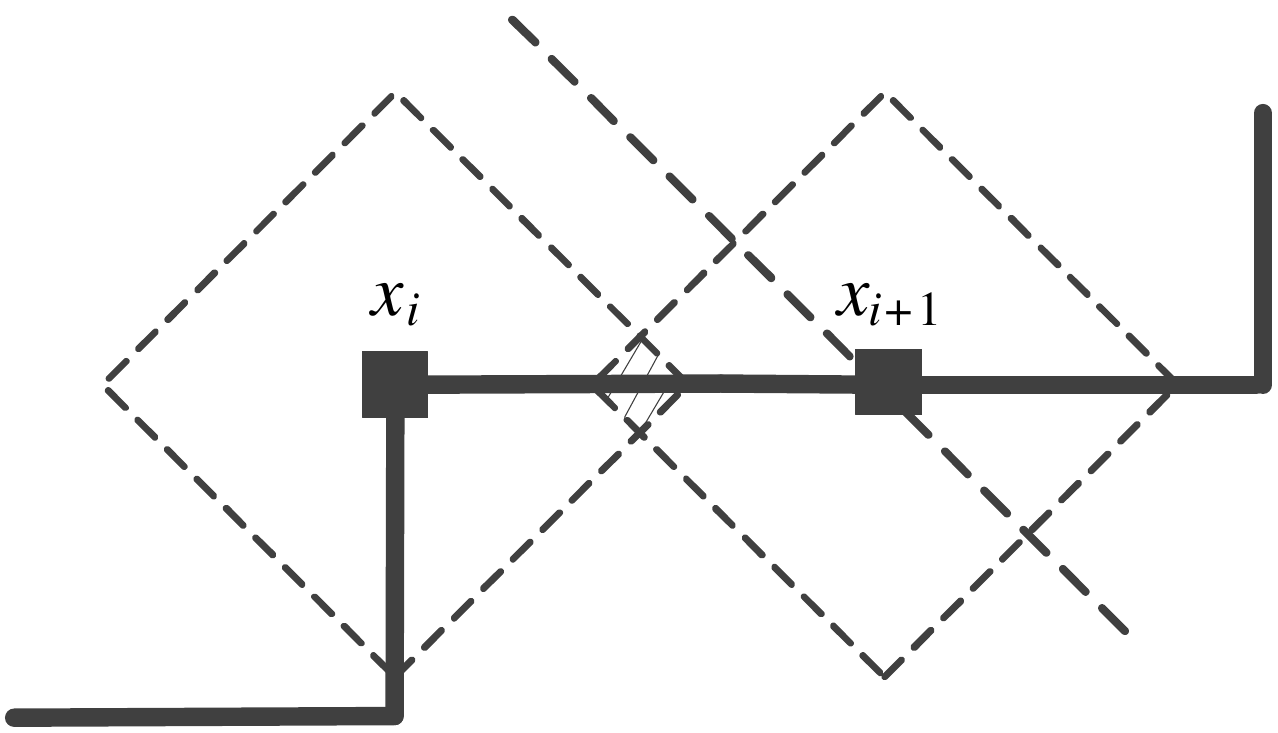}}
 \hspace{0.1in}
 \subfigure[]{
 \label{Figure_Line_Square} 
 \includegraphics[width=0.3\textwidth]{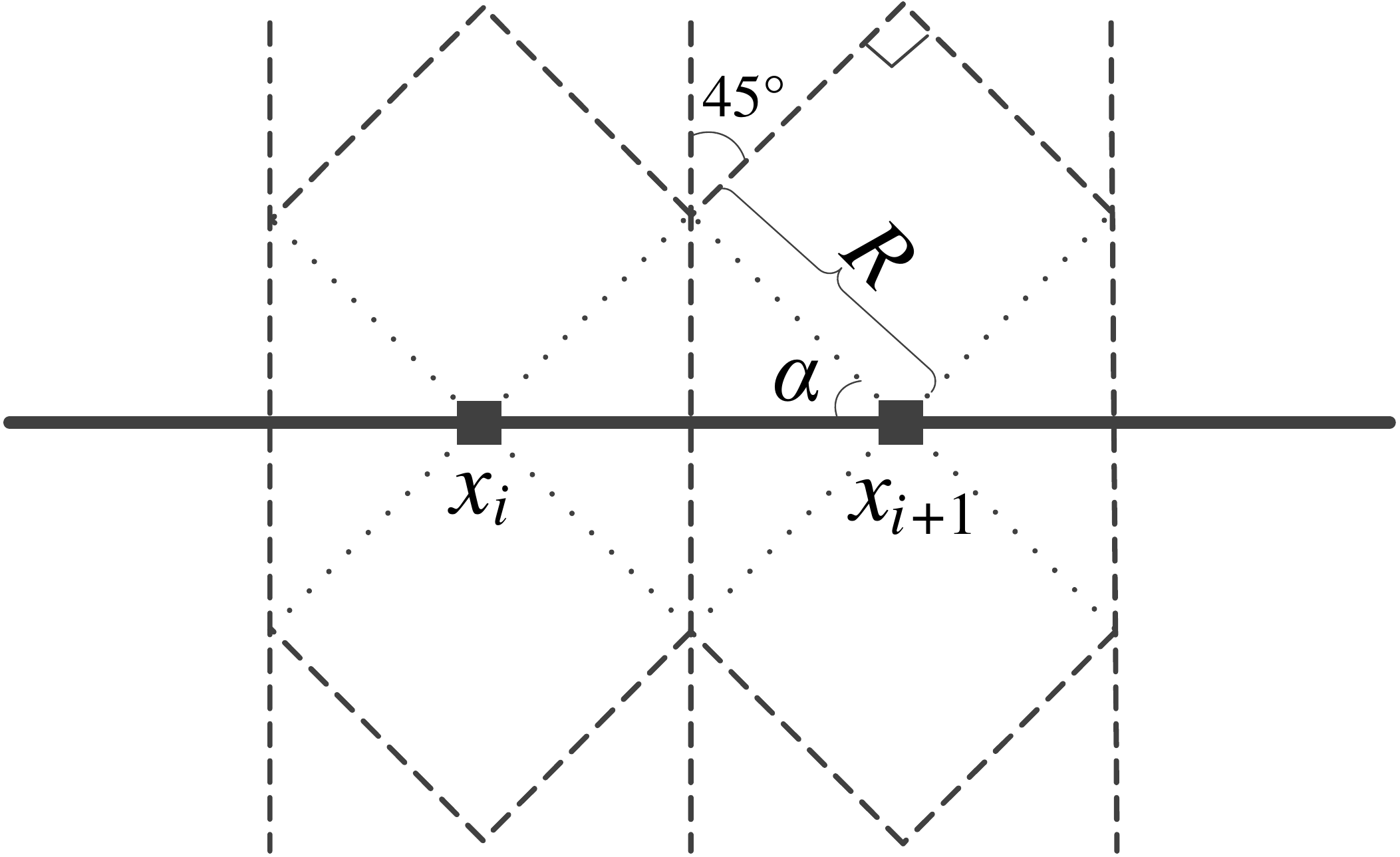}}
 \caption{Construction of a lower bound under $L_1$ metric}
 \vspace*{-10pt}
 \label{Fig_Lem_Given_L_N_A_Boundary_L_1} 
\end{figure}

Similar to Section \ref{sec_LB_L_2}, we now construct a cost lower bound by considering all solutions that incur a fixed inbound travel length $l$, a fixed number of facilities $N$ that collectively cover the customers in an area of total size $A$. Note from Section~\ref{L_1_C_0} that the square shape minimizes the outbound cost for any given size of service regions. Thus, if the inbound travel length $l$ is larger than the total diagonal length of $N$ identical squares (each with area size $\frac{A}{N}$), then this case degrades to a trivial one where the optimal cost is achieved when all $N$ service regions take the shape of identical squares, as shown in Figure~\ref{Figure_Zigzag_L_1}. However, as discussed before, this case never yields a good cost lower bound since we can always shift the facility locations and their service regions along the TSP tour to reduce the length of inbound truck.

Now we consider the more general case when $l\leq \left(2NA\right)^{\frac{1}{2}}$. Overlaps among neighboring service regions are now inevitable, forcing facilities to serve farther customers (outside of the ideal squares of size $\frac{A}{N}$) and incur higher outbound cost. As such, minimizing the total overlapping area is equivalent to minimizing the total outbound cost for any given $l, N, A$. Note that for any two overlapping service areas, say $i$ and $i+1$ as shown in Figure~\ref{Figure_Zigzag_Overlap_L_1}, the overlap area will be minimized if we shift facility $i+1$ (and all facilities $i+2, i+3, \cdots$) along the 45 degree line so that segment $x_i x_{i+1}$ becomes parallel to one of the coordinate axes as shown in Figure~\ref{Figure_Zigzag_Overlap_Line_L_1}. The inbound TSP distance will not change, but the outbound cost will decrease. We can repeat this for all neighboring facilities along the TSP tour, and a lower outbound cost for any given $l, N, A$ will be achieved when all $N$ facilities are along a straight line parallel to a coordinate axis.

The rest of the argument is very similar to that in Section \ref{sec_LB_L_2}. When all facilities are along a straight TSP line, any two adjacent service regions should be separated by a boundary line that is perpendicular to the TSP tour (see Figure~\ref{Figure_Line_Square}); each facility will only serve the customers within the two nearest boundaries. To minimize the outbound cost of each facility within its boundaries, the optimal service region should be the intersection of the area between the two boundaries and a square shape; otherwise we can always perturb customer allocation to reduce outbound cost. Also note that the optimal locations and service regions of the facilities should form a centroidal Voronoi tessellation, as shown in Figure~\ref{Figure_Line_Square}; i.e., any facility should also be at the center of its service region. Hence, all service regions must be identical, and all facilities are evenly spaced along the straight TSP tour.


The above lower bound (mainly regarding outbound cost) holds for a fixed inbound cost (i.e., given values of $l, N, A$). The best lower bound can be obtained by choosing proper values of $l, N, A$ to address the inbound and outbound cost trade-off. We consider the geometry in Figure~\ref{Figure_Line_Square} and express $l$ as a function of $N,A$ and $\alpha$; i.e.,
\[
l=\left(2AN\right)^{\frac{1}{2}}\left(1+2\tan\alpha\right)^{-\frac{1}{2}}.
 \]
Following the notation in Figure~\ref{Figure_Line_Square}, the best lower bound can be achieved by finding the optimal $R$ and $\alpha$ that solves the following problem:
\begin{equation}
\min \frac{Nf}{A}+\frac{\kappa f N}{A}\left(\frac{2 }{3}R^3\cos\alpha\left(2+\sin^2\alpha+6\sin\alpha\cos\alpha\right)+2r R\frac{A}{N}\cos\alpha\right),\label{eq_lower_bound_L_1}
\end{equation}
subject to $2NR^2\left(\cos^2\alpha+\sin2\alpha\right)=A$. 

It is easy to show from the first order condition of \eqref{eq_lower_bound_L_1} with respect to $\alpha$, that there is a single optimizer $\alpha^*=\arctan (2r+\sqrt{(2r+4r^2)}) \in [0, \frac{\pi}{2})$. We write it as a function of $r$, and define a new function
\begin{equation}
\bar{g}(r )=\frac{3\left(2\cos\alpha^*(r)\right)^\frac{1}{2}\left(2\sin\alpha^*(r) +\cos\alpha^*(r)\right)^{\frac{3}{2}}}{3\sin \left(2\alpha^*(r)\right)-2\cos\left(2\alpha^*(r)\right)+4}. \label{bar_g_r_inbound_L_1}
\end{equation}
Then we have
\allowdisplaybreaks\begin{align*}
\eqref{eq_lower_bound_L_1}\geq \frac{Nf}{A}+\frac{\kappa f A^{\frac{1}{2}}}{\sqrt{N}\bar{g}\left(r \right)}
\geq 3\sqrt[3]{\frac{\kappa^2f^3}{4\bar{g}^2(r )}}.
\end{align*}
The last inequality becomes equality by choosing $\frac{A}{N}=\left(\frac{\kappa}{2}\right)^{-\frac{2}{3}}\left(\bar{g}(r )\right)^{\frac{2}{3}}$.

We further note that the elongated hexagons in Figure~\ref{Figure_Line_Square}, which achieve the cost lower bound, can also form a feasible spatial tessellation. Hence, it yields an optimal tessellation. This finding is summarized in the following proposition.

\begin{proposition}\label{P_LB_inbound_L_1}
Under $L_1$ metric, the optimal shape of facility service region is elongated hexagon. With $\alpha^*(r)=\arctan (2r+\sqrt{(2r+4r^2)})$ and $\bar{g}(r )$ as defined in \eqref{bar_g_r_inbound_L_1}, the optimal value of \eqref{original_obj} is given by
$3\sqrt[3]{\frac{\kappa^2f^3}{4\bar{g}^2(r )}}$.
\end{proposition}

\section{Discussion}\label{sec_Discussion}

\subsection{Sensitivity Analyses}
\begin{figure}[ht]
 \centering
 \subfigure[$\alpha^*(6,r ),\bar{\alpha}^*(6,r )$ and $\alpha^*(r )$ versus $r $]{
 \label{t_optimal} 
 \includegraphics[width=0.45\textwidth]{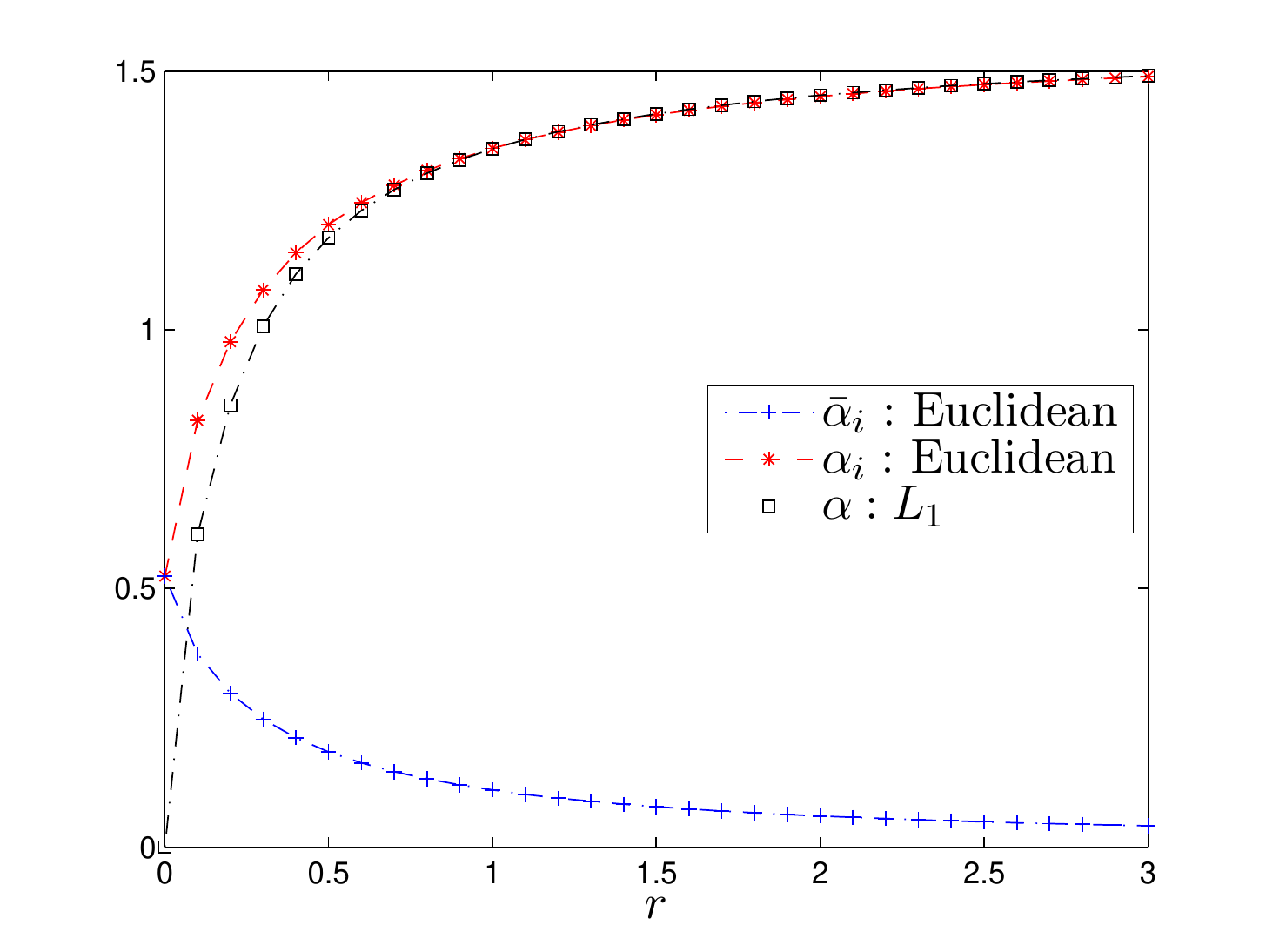}}
 \hspace{0.1in}
 \subfigure[$z_{ub}^*$ and $z^*$ versus $r $]{
 \label{value_of_Optimal}
 \includegraphics[width=0.45\textwidth]{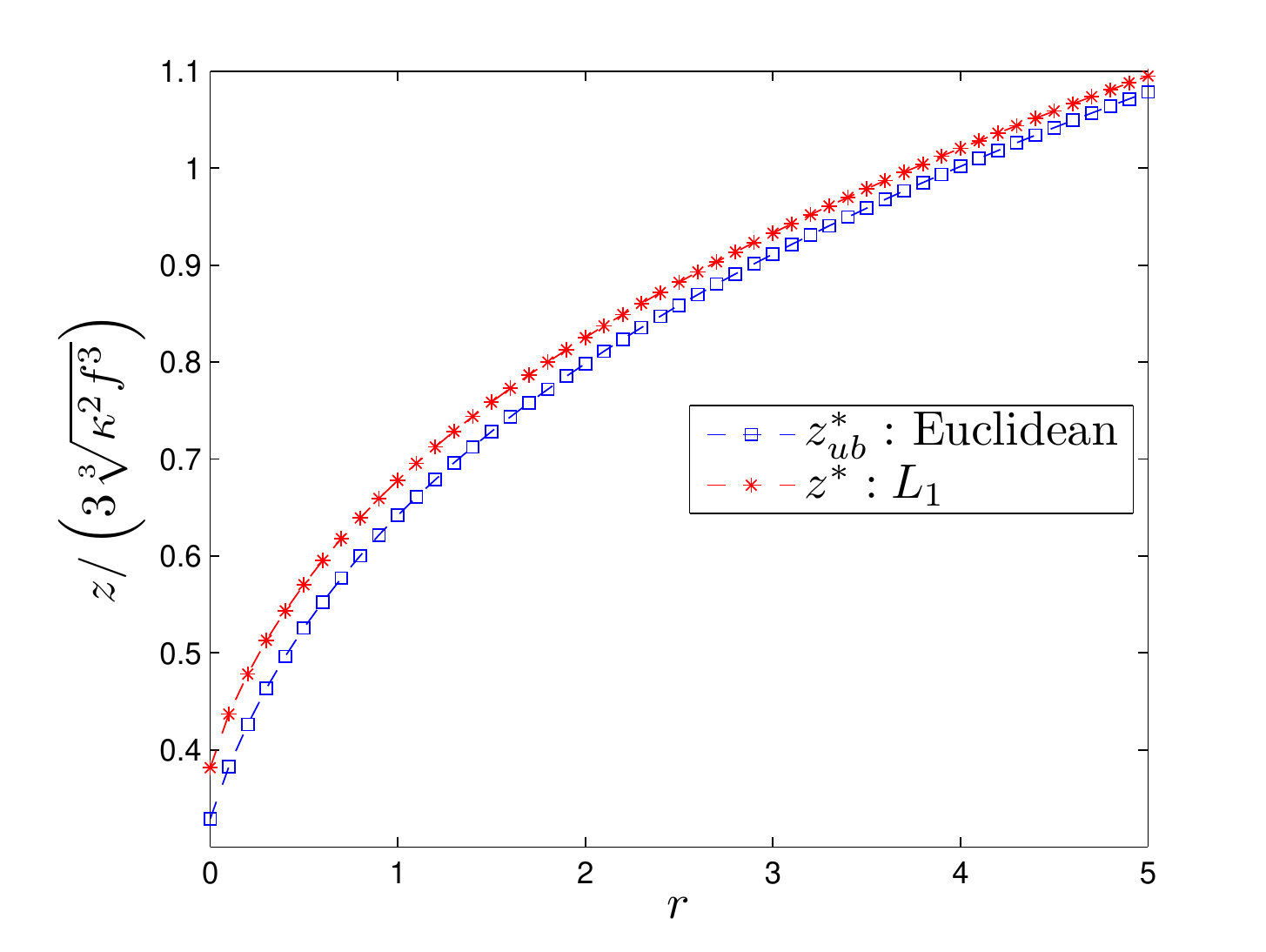}}
 \caption{Impact of cost coefficient ratio $r$}
 \label{verify_optimal_in} 
\end{figure}

Figure~\ref{t_optimal} plots the optimal basic angles (i.e., $\alpha_i,\bar{\alpha}_i$ for Euclidean metric, $\alpha$ for $L_1$ metric) as functions of $r $. Figure \ref{value_of_Optimal} plots the upper bound for Euclidean metric, $z^*_{ub}$, and the optimal solution for $L_1$ metric, $z^*$, as functions of $r $. We notice that all functions are monotone, and for sufficiently large $r $, $\alpha_i \rightarrow \frac{\pi}{2},\bar{\alpha}_i\rightarrow 0, \alpha\rightarrow \frac{\pi}{2}$ and $z^*_{ub}, z^*$ increase concavely. This means when the inbound transportation cost is very high, the facility service shape will be elongated to shorten the inbound truck travel. This finding is consistent with, but generalizes, the asymptotic results in \cite{Carlsson13}. We also note that $z^*>z_{ub}^*$ for all $r \ge 0$, which is intuitive because the Euclidean distance for any two points in a plane is no larger than their $L_1$ distance. However, the cost difference diminishes when $r$ becomes larger (i.e., when inbound cost dominates), because the optimal service regions become very thin stripes (regardless of the metric), and hence the influence of the distance metric becomes insignificant.

\subsection{Numerical Verification}\label{NS}

To further verify our cost bounds, we solve a discrete mathematical program using a grid of $M\times M$ points denoted by set $\mathbb{G}$ ($\left\vert\mathbb{G}\right\vert=M^2$). Each point $i\in \mathbb{G}$ denotes a customer as well as a candidate facility location. An arbitrary point $o\in\mathbb{G}$ represents the start point of the TSP tour. Here, we use $X_i=1$ to indicate that the facility will be built at point $i\in \mathbb{G}$ with cost $f_i$; otherwise, $X_i=0$. In terms of outbound delivery, we use $Y_{i,j}=1$ to denote when customer $j\in\mathbb{G}$ will be assigned to facility $i\in\mathbb{G}$ with cost $c_{ij}d_{i,j}$ ($d_{i,j}$ denotes the Euclidean or $L_1$ distance between $i$ and $j$). Meanwhile, as for the TSP tour, $Z_{i,j}=1$ denote facility $j\in\mathbb{G}$ is visited right after facility $i\in\mathbb{G}$ by the inbound truck at cost $C_{i,j}d_{i,j}$. We use additional continuous variable $u_i, i\in\mathbb{G}$, to avoid sub-tours.

As such, the following mixed-integer program represents the discrete version of our problem.
\allowdisplaybreaks\begin{subequations}\begin{align}
\min\limits_{\bf{X,Y,Z}}&\sum\limits_{i\in \mathbb{G}}f_iX_i+\sum\limits_{i,j\in \mathbb{G}}\left(c_{ij}d_{i,j}Y_{i,j}+C_{ij}d_{i,j}Z_{i,j}\right)\label{eq:obj_IP}\\
\text{s.t.}\,\,\,\,
&\sum\limits_{i\in \mathbb{G}}Y_{i,j}=1,\forall j\in \mathbb{G}\label{con:Y_sum_1_IP}\\
&Y_{i,j}\leq X_i,\forall \left(i,j\right)\in \mathbb{G}\times\mathbb{G}\label{con:y<x_IP}\\
&\sum\limits_{i\in \mathbb{G}}Y_{i,j}=1,\forall j\in \mathbb{G}\label{con:y=1_IP}\\
&\sum\limits_{i\in \mathbb{G}}Z_{i,j}=X_i,\forall j\in \mathbb{G}\label{con:z_x_1_IP}\\
&\sum\limits_{j\in \mathbb{G}}Z_{i,j}=X_j,\forall j\in \mathbb{G}\label{con:z_x_2_IP}\\
&u_i-u_j+M^2Z_{i,j}\leq M^2-1,\forall \left(i,j\right)\in \mathbb{G}\times\mathbb{G}\backslash\left\{\left(0,0\right)\right\}\label{con:overcome_subtour_IP}\\
&0\leq u_i\leq M,\forall i\in \mathbb{G}\backslash\left\{0\right\};u_o=0\label{con:u_range_IP}\\
&X_i\in\left\{0,1\right\},\forall i\in\mathbb{G};Y_{i,j},Z_{i,j}\in\left\{0,1\right\},\forall \left(i,j\right)\in \mathbb{G}\times\mathbb{G}.\label{con:binary_variable_IP}
\end{align}
\end{subequations}
Here, objective function~\eqref{eq:obj_IP} minimizes the total facility set-up cost, outbound delivery cost and inbound transportation cost. Constraints~\eqref{con:y<x_IP} and~\eqref{con:y=1_IP} postulate that each customer should be sent to one built facility. Constraints~\eqref{con:z_x_1_IP} and~\eqref{con:z_x_2_IP} ensure that each facility should be passed by the inbound truck. Constraints~\eqref{con:overcome_subtour_IP} eliminate sub-tours. Constraints~\eqref{con:u_range_IP} and~\eqref{con:binary_variable_IP} define the continuous and binary variables.

\begin{figure}[htbp]
 \centering
 \subfigure[Euclidean metric]{
 \label{Numerical_result} 
 \includegraphics[width=0.45\textwidth]{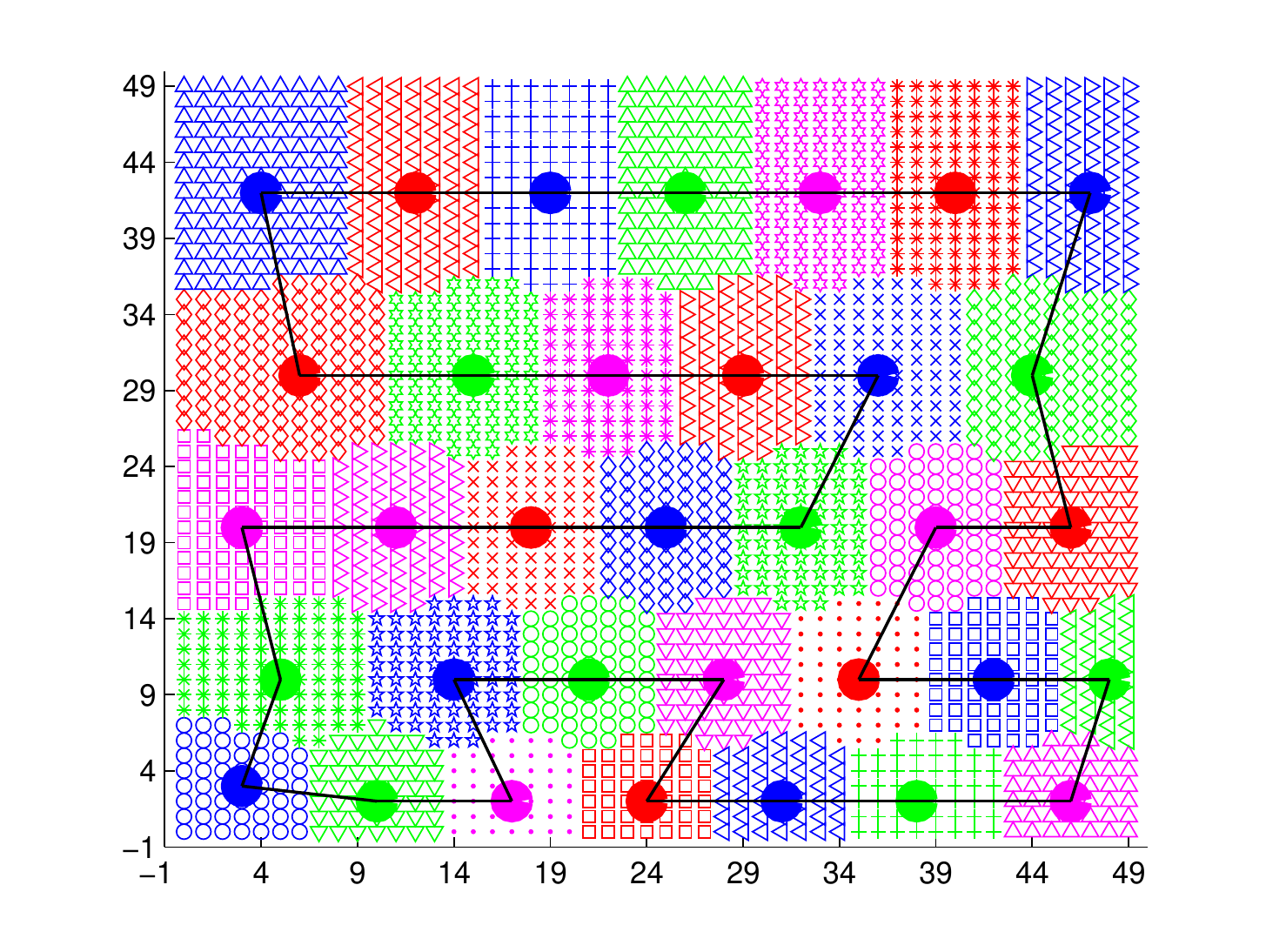}}
 \hspace{0.1in}
 \subfigure[$L_1$ metric]{
 \label{Numerical_result_L_1}
 \includegraphics[width=0.45\textwidth]{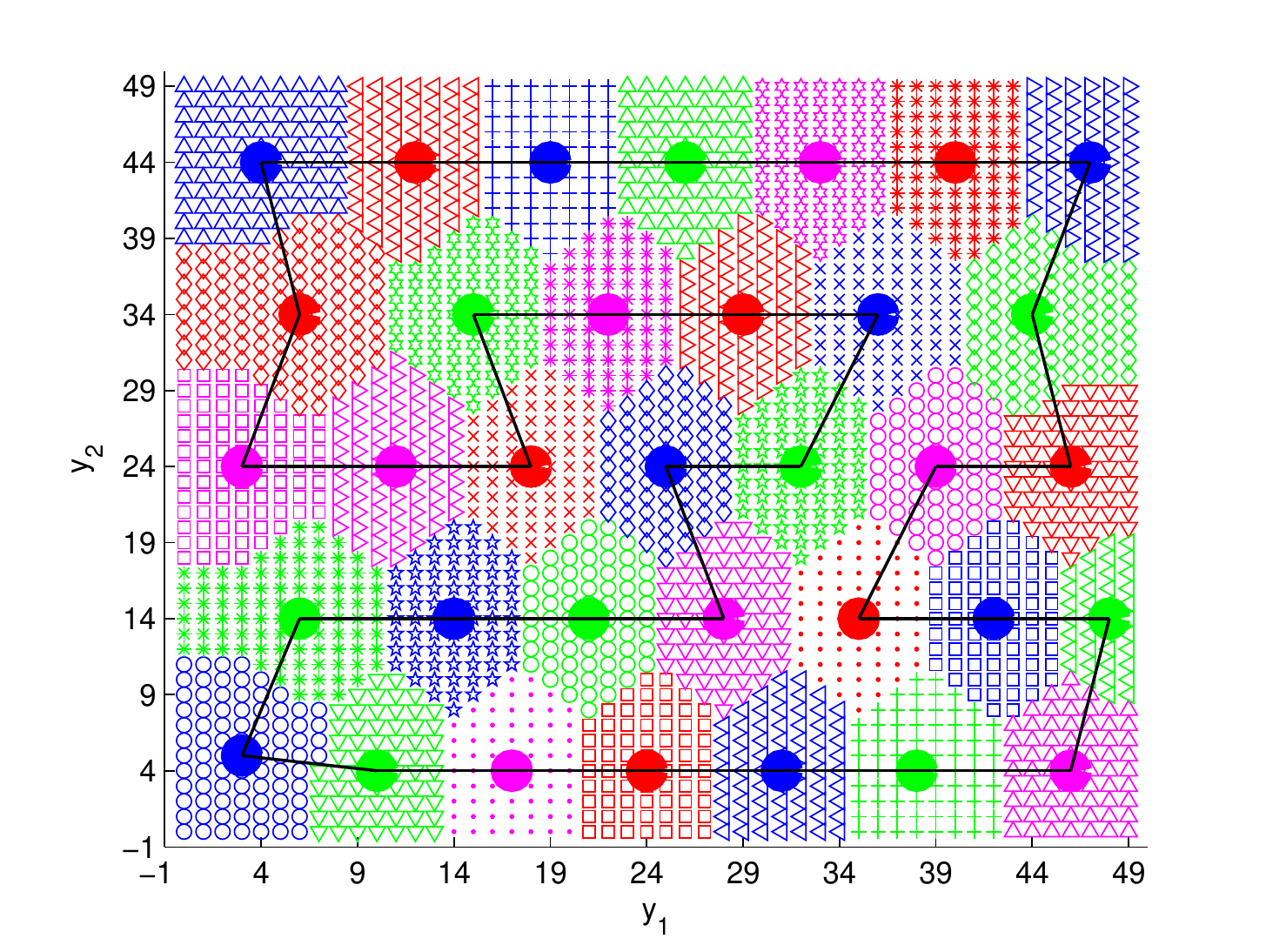}}
\caption{Results of numerical experiments on a discrete grid}
 \vspace*{-5pt}
 \label{Numerical_Experiment} 
\end{figure}

Since the discrete problem is very difficult, we use a simulated annealing heuristic to solve it. Figure~\ref{Numerical_result} shows the computation result for an instance under Euclidean metric, with $M=50$, $f_i=299.66,c_{ij}=1,C_{ij}=12,\forall i,j\in\mathbb{G}$. We notice that most of facility service region shapes (especially those away from the boundaries) turn out to be cyclic hexagonal, exactly as what we would expect from Proposition~\ref{Prop_lB_Z_N_A_inbound}. Figure~\ref{Numerical_result_L_1} shows the result under $L_1$ distance metric, with the same parameters except for $f_i=199.31,\forall i\in\mathbb{G}$ and $L_1$ distances. The facility service regions turn out to have noncyclic hexagonal shapes, again, as expected. For both cases, the finite boundaries of the $50 \times 50$ area do seem to influence some of the service region shapes, but such effect shall diminish when $M \rightarrow \infty$.

We also measure from Figure~\ref{Numerical_result} the average basic angles for those cyclic hexagons in the center of the region, which turn out to be $18.8^\circ$ and
$52.3^\circ$. From \eqref{equation_alpha_1_2_V2_H} and \eqref{total_cost_system_opt_n_in}, the theoretical number of facilities $N^*\approx 34$, and the optimal basic angle values are $18.4^\circ$ and $53.2^\circ$.
We can see that the error between the theoretical result and the experimental measurement is only $0.9^\circ$. Similarly, the average basic angles for the non-cyclic hexagons in Figure~\ref{Numerical_result_L_1} are around $45.0^\circ$, while the theoretical value is $44.5^\circ$, yielding a $0.5^\circ$ difference. We anticipate that these minor errors will further diminish if we increase the size of the grid (so as to eliminate the influence of the boundaries), and discretize it into a finer resolution.

\subsection{Impacts of Inventory Cost}
Inventory cost is sometimes significant to a transshipment system. Attempts have been made to incorporate inventory considerations into discrete facility location models~\citep{Daskin02,Chen2011joint} and location-routing models~\citep{Shen2007incorporating}. In this section, we will consider cost for inbound inventory holding in transshipment facilities and discuss its impacts on the optimal system design \eqref{original_obj}. Since the demand rate for each facility is deterministic, constant over time, and proportional to the size of the service region, we can adopt a simple cycle inventory policy to determine the optimal inbound replenishment frequency, and assume that the fixed order cost (per order) and inventory holding cost (per item-time) are both constants. To simplify the formulas, however, we express these cost coefficients in the form of $b f\kappa^{\frac{1}{3}}/\lambda$ and $h f\kappa^{\frac{1}{3}}$, respectively, for some proper constants $b$ and $h$. As such, $b$ and $h$ essentially represent the relative magnitudes of the fixed order cost and inventory holding cost as compared to the other costs (e.g., facility set-up cost, transportation cost). EOQ (economic order quantity) trade-off can be directly applied to determine replenishment frequency, and the total inventory cost of the $i$th facility is $f\left(2bhA_i\right)^{\frac{1}{2}}\kappa^{\frac{1}{3}}$. Thus, \eqref{original_obj} is now generalized as follows:
\allowdisplaybreaks\begin{align}
z = &\min\limits_{{\cal{N}},\left\{x_i\right\},\left\{{\cal{A}}_i\right\}} \lim_{\begin{subarray}{l} N\rightarrow\infty\end{subarray}}
{\frac{f}{\sum\limits_{i=1}^{N} {A_{i}}}\sum\limits_{i=1}^{N} {\left(1+\kappa\int_{{\cal{A}}_{i}}{ \|x-x_i\| d{x}}+\left(2bhA_i\right)^{\frac{1}{2}}\kappa^{\frac{1}{3}}\right)}+\frac {\kappa r}{N} \sum\limits_{i=1}^{N} {l_i}} \label{original_obj_inventory}.\end{align}

We shall first note that the inventory cost term $f\left(2bhA_i\right)^{\frac{1}{2}}\kappa^{\frac{1}{3}}$ is concave with respect to $A_i$. The optimal spatial tessellation patterns discussed in previous sections may no longer hold if the value of $bh$ is large; e.g., in the extreme, if inventory cost dominates outbound and inbound transportation cost, it is beneficial to aggregate demand, and hence the optimal design should degenerate to a single facility that serves all customers on the plane. Hence, in general, inventory cost could have a significant impact on the optimal spatial facility layout.

However, we shall also note from the proofs of all lemmas and theorems in Sections \ref{Without inbound truck} - \ref{sec_configuration_L_1} that the optimality of the hexagonal spatial tessellations and facility layouts (under both metrics) will still hold as long as the objective function \eqref{original_obj_inventory} remains convex with respect to $A_i, \forall i$. Such a condition could be achieved in many ways. First, it suffices if the parameter $bh $ becomes facility specific, i.e., $(bh)_i$, and it is dependent on $A_i$ and satisfies the following second-order condition
\[
2(bh)_iA_i^2\frac{d^2\left(bh\right)_i}{dA_i^2}-\left((bh)_i-A_i\frac{d\left(bh\right)_i}{dA_i}\right)^2\geq 0.
\]
For example, the fix order cost coefficient could be a convex function of the total served demand~\citep{Veinott1964production}, while the holding cost coefficient remain constant; e.g., for the $i$th facility, $(bh)_i=A_i^\beta$, for some $\beta\geq 1$. The inventory cost term in \eqref{original_obj_inventory} becomes $\sqrt{2}f\kappa^{\frac{1}{3}}A_i^{\frac{1+\beta}{2}}$, which is now convex in $A_i$.

Another obvious condition for \eqref{original_obj_inventory} to remain convex is that the inventory cost term be dominated by other costs. In this case, the best tessellation for Euclidean metric and $L_1$ metric remain elongated hexagons (as shown in previous sections), but the optimal value of $A/N$ may change. Simple algebra shows that closed-form formulas can still be found as follows:
\allowdisplaybreaks\begin{align}
\frac{A}{N}=\left\{\begin{array}{cl}
\frac{4\sqrt{2}}{3}\kappa^{-\frac{2}{3}}g\left(n,r \right)\left(bh\right)^{\frac{1}{2}}\cos^2\left[\frac{1}{3}\arccos\left(\left(2bh/9\right)^{-\frac{3}{4}}g^{-\frac{1}{2}}\left(n,r \right)\right)\right],n=6,\infty,&\text{Euclidean metric},\\
\frac{4\sqrt{2}}{3}\kappa^{-{\frac{3}{2}}}\bar{g}\left(r \right)\left(bh\right)^{\frac{1}{2}}\cos^2\left[\frac{1}{3}\arccos\left(\left(2bh/9\right)^{-\frac{3}{4}}\bar{g}^{-\frac{1}{2}}\left(r \right)\right)\right],&L_1\text{ metric}.
\end{array}\right. \label{A_N_inventory}\end{align}

\begin{figure}[htbp]
 \centering
 \subfigure[Euclidean metric]{
 \label{Pencent_Dif_Z} 
 \includegraphics[width=0.45\textwidth]{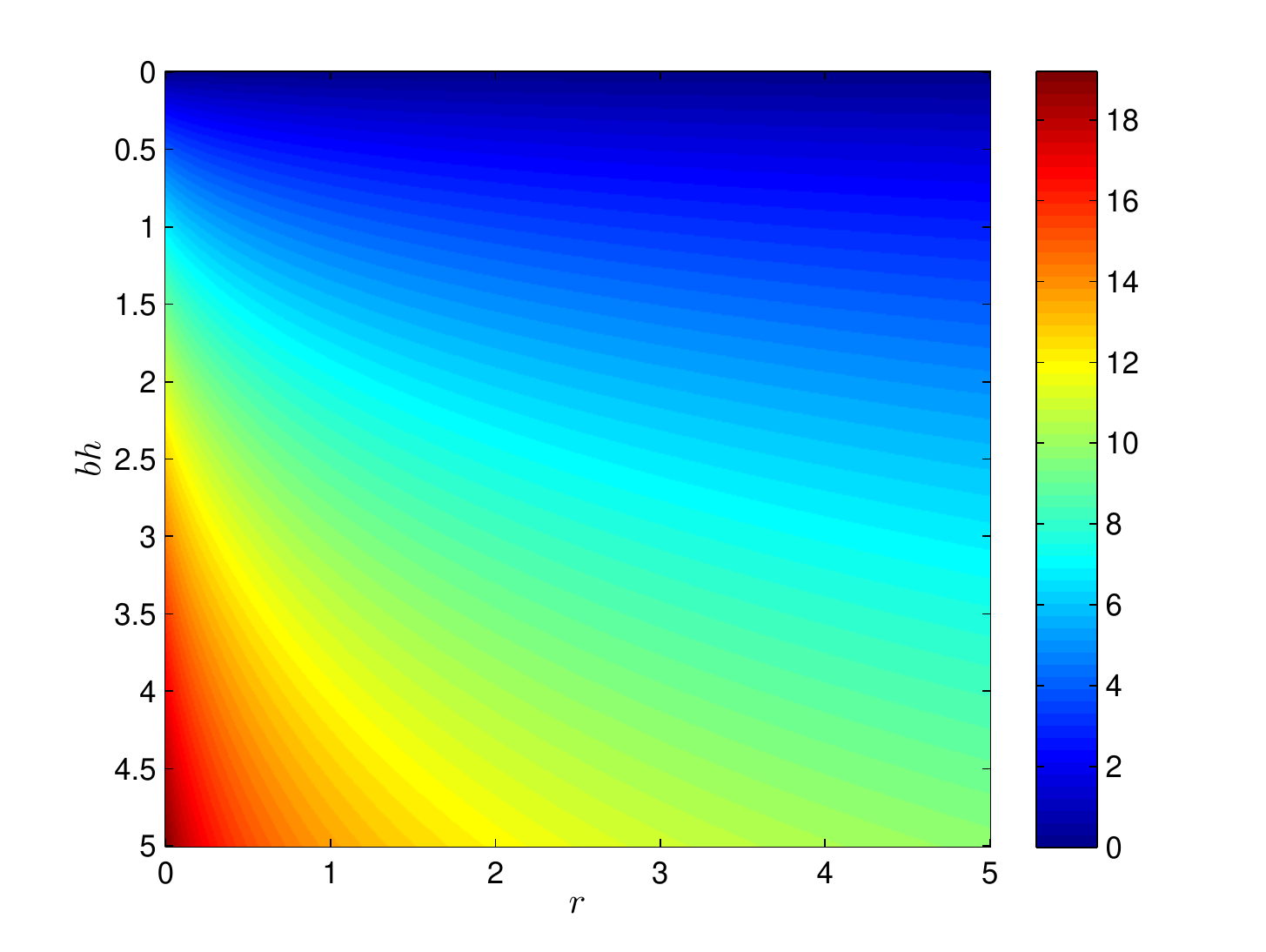}}
 \hspace{0.1in}
 \subfigure[$L_1$ metric]{
 \label{Pencent_Dif_Z_L_1}
 \includegraphics[width=0.45\textwidth]{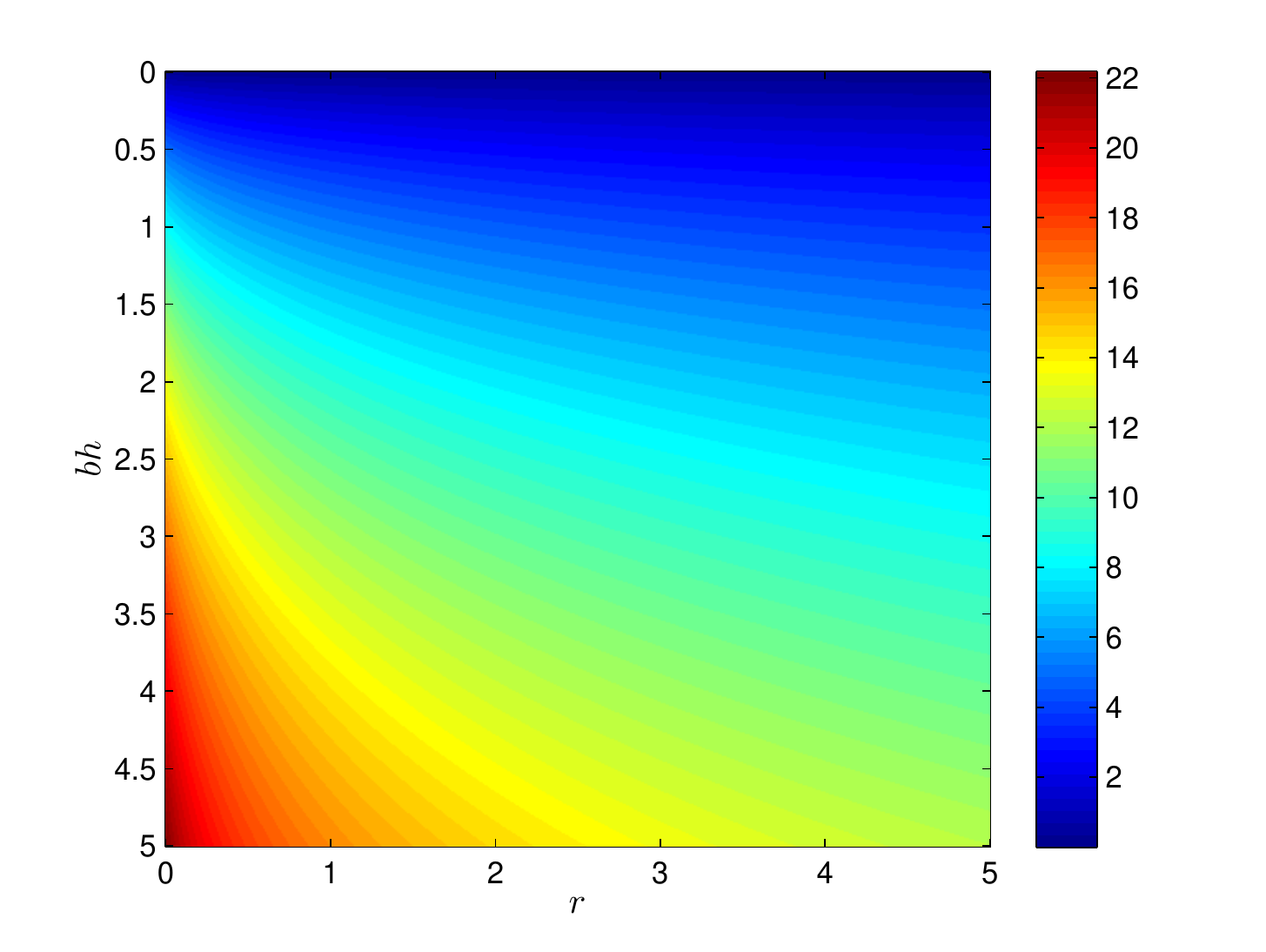}}
\caption{{Percentage difference in system cost with and without inventory cost consideration.}}
 \vspace*{-5pt}
 \label{Pencent_Dif_Z_Inventory} 
\end{figure}

Figure~\ref{Pencent_Dif_Z_Inventory} illustrates the percentage differences in objective function \eqref{original_obj} between the optimal solutions with and without considering inventory cost. It turns out that inventory cost has a slightly larger effect under $L_1$ metric than under Euclidean metric. We also note that for both metrics, when $r$ grows, the inventory cost becomes less significant, so the cost difference becomes smaller. Conversely, when the value of $bh$ grows larger, inventory cost plays a more important role and leads to a bigger cost difference.

\section{Conclusion}\label{conclusion}
This paper studies the optimal transshipment facility layout on a homogeneous plane $\Re^2$ that minimizes total system cost for facility set-up, outbound customer delivery and inbound replenishment transportation. We first show a proof for Gersho's conjecture~\citeyearpar{Gersho79} under Euclidean metric, which states that when inbound transportation cost is negligible, the optimal spatial partition of $\Re^2$ should be regular hexagons. When inbound cost is non-negligible, we first construct an upper bound by tessellating the plane with elongated cyclic hexagons. Then we derive a cost lower bound which is achieved by an infeasible infinite cyclic polygon. We derive analytical formulas for both bounds, and show that the percentage gap between these two bounds is quite small (i.e., within 0.3\%). We further illustrate the impact of service region shapes by comparing the performance of three special shapes (i.e., triangle, rectangle, hexagon), and show that elongated cyclic hexagon outperforms the others. To verify our analytical results, we also formulate a mixed-integer program locating-routing model and observe the near-optimal spatial tessellation via a numerical experiment on a grid of points. The numerical results turn out to be very close to our analytical predictions. Finally, we extend our discussion to the $L_1$ metric case and show that a similar non-cyclic hexagon shape becomes exactly optimal when they are properly oriented along the axes of the $L_1$ coordinate system.

Future research can be conducted in several directions. For example, we may introduce a finite capacity of the inbound truck and investigate its impact on the optimal spatial layout. Moreover, we strongly suspect that the elongated cyclic hexagonal shape is actually optimal under the Euclidean metric, but we have only presented it as a feasible solution (which yields an upper bound). It will be ideal to further prove the optimal tessellation pattern for the Euclidean metric and other variations (e.g., other metrics). Finally, we have shown that when inventory cost becomes dominant, the spatial tessellation patterns presented in this paper may no longer be optimal. It will be interesting to find the optimal tessellation pattern for more general settings.

\section*{Acknowledgments}
This research was supported in part by the U.S. National Science Foundation through Awards EFRI - RESIN - 0835982, CMMI - 1234085 and CMMI - 0748067. The helpful comments from Professors Carlos Daganzo, Max Zuo-Jun Shen (UC Berkeley) and James Campbell (U of Missouri at St. Louis) on an earlier version of the paper are gratefully acknowledged.

\bibliography{Facility}
\newpage

\begin{appendices}
\addcontentsline{toc}{chapter}{Appendices}
\renewcommand\thesection{Appendix A}
\section{Proofs of the Lemmas and Propositions}


\subsection{Proof of Lemma~\ref{L1_Given_area_angle}}\label{proof_L1_given_area_angle}
\begin{figure}[htbp]
 \begin{center}
 \includegraphics[width=0.3\textwidth]{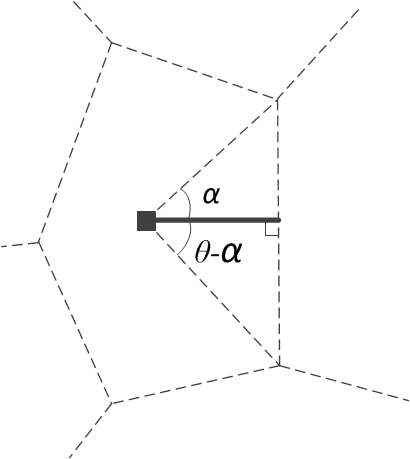}
 \caption{Notation for the proof of Lemma~\ref{L1_Given_area_angle}} \label{figure1}
 \vspace*{-20pt}
 \end{center}
\end{figure}

\begin{proof}
Given the area of a basic triangle $A$ and its basic angle $\theta$, without losing generality, we assume that $\alpha\leq{\frac{\theta}{2}}<\frac{\pi}{2}$, as shown in Figure~\ref{figure1}. Simple algebra will show that the outbound delivery cost in this triangle can be formulated as a function of $\alpha$, i.e.,
\begin{equation}
 z_{\Delta}\left(\alpha\right)
 =\frac{\kappa f }{3}\left(\sin\theta\right)^{-\frac{3}{2}} A^\frac{3}{2}\left(\cos\theta+\cos\left(\theta-2\alpha\right)\right)^\frac{3}{2}
 \int_{-(\theta-\alpha)}^{\alpha}\frac{1}{\cos^3 t} dt,\label{z_v_d}
\end{equation}
and $\frac{d z_{\Delta}\left(\alpha\right)}{d\alpha}<{0}$ when $0<\alpha<{\frac{\theta}{2}}$. This shows that $\alpha=\frac{\theta}{2}$ is the unique solution to minimize $z_{\Delta}\left(\alpha\right)$. This completes the proof.
\end{proof}

\subsection{Proof of Lemma~\ref{P_Given_Area_All}}\label{proof_P_Given_Area_All}
\begin{proof}
Let $F(x)=\kappa f\cos^3x\int_{-x}^{x}\frac{1}{\cos^3 t} dt$. In order to minimize the total outbound delivery cost, we move $\sum\limits_{j=1}^nR_j^2\sin\alpha_j\cos\alpha_j=A$ into~\eqref{z_all} with Lagrangian multiplier $\gamma$:
\begin{align}
 \sum\limits_{j=1}^n \frac{R_j^3}{3}F(\alpha_j)+\gamma(A-\sum\limits_{j=1}^nR_j^2\sin\alpha_j\cos\alpha_j).\label{z_all_Lag}
\end{align}
First order condition shows that
\begin{align}
 R_j^2F(\alpha_j)-2\gamma R_j\sin\alpha_j\cos\alpha_j =0, \forall j.\label{z_all_Lag_De}
\end{align}
which yields
\begin{equation*}
R_j=\frac{2\gamma\sin\alpha_j\cos\alpha_j}{F(\alpha_j)},\forall j \text{ and } \gamma =\frac{A^{\frac{1}{2}}}{2}\left(\sum\limits_{j=1}^n\frac{\sin^3\alpha_j\cos^3\alpha_j}{F^2(\alpha_j)}\right)^{-\frac{1}{2}}.\label{s_R}
\end{equation*}

We can further show that
\begin{equation}
 \eqref{z_all}\geq \frac{1}{3}\kappa f A\sqrt{A}\left(\sum\limits_{j=1}^n\frac{\tan^3\alpha_j}{\left(\int_{-\alpha_j}^{\alpha_j}\frac{1}{\cos^3 t} dt\right)^2}\right)^{-\frac{1}{2}}.\label{z_all_No_R}
\end{equation}

Consider the following function:
\begin{equation}
 \phi(x)= \frac{\tan^3x}{\left(\int_{-x}^{x}\frac{1}{\cos^3 t} dt\right)^2}
= \frac{\tan^3x}{\left(\log\tan\left(\frac{x}{2}+\frac{\pi }{4}\right)+\frac{\tan x}{\cos x}\right)^2},\label{func_h}
\end{equation}
By checking the first order and second order derivatives of $\phi(x)$, we can find $\phi(x)$ to be strictly concave over $x\in \left(0,\frac{\theta }{2}\right)$. Therefore, Jensen's inequality yields
\begin{align}
\sum\limits_{j=1}^n\frac{\tan^3\alpha_j}{\left(\int_{-\alpha_j}^{\alpha_j}\frac{1}{\cos^3 t} dt\right)^2}=n\sum\limits_{j=1}^n\frac{\phi(\alpha_j)}{n}\leq{n\phi\left(\sum\limits_{j=1}^n\frac{\alpha_j}{n}\right)}=n\phi\left(\frac{\theta }{n}\right).\label{concave_h}
\end{align}

Substituting~\eqref{concave_h} into~\eqref{z_all_No_R}, we can get the following inequality:
\begin{equation}
 \eqref{z_all}\geq{ \frac{\kappa f A\sqrt{A}}{3\sqrt{n\phi\left(\frac{\theta }{n}\right)}}}= \frac{1}{3}\kappa f A\sqrt{A}\left(\log\tan\left(\frac{\theta }{2n}+\frac{\pi }{4}\right)+\frac{\tan \frac{\theta }{n}}{\cos \frac{\theta }{n}}\right)\left(n\tan^3\frac{\theta }{n}\right)^{-\frac{1}{2}}.\label{lower_z_all}
\end{equation}

The equality holds only if $\alpha_j=\frac{\theta }{n}$ and $R_j=\sqrt{\frac{2A}{n\sin\frac{2\theta }{n}}},\forall j$.
\end{proof}

\subsection{Proof of Lemma~\ref{Lem_lB_Z_N_A}}\label{proof_Lem_lB_Z_N_A}
\begin{proof}
First of all, since $g^2(x)$ is strictly concave and monotonically increasing over $\left[3,+\infty\right)$, Jensen's inequality leads to the following:
\begin{equation}
\sum\limits_{i=1}^{N} {g^2(n_i)}\leq{Ng^2\left(\sum_{i=1}^N n_i/N\right)}.\label{opt_lb}
\end{equation}
We now try to minimize the right hand side of \eqref{total_cost_system} as an unconstrained optimization problem by adding constraint~\eqref{Area_con_system} into the objective with Lagrangian multiplier $\mu$; i.e.,
\begin{equation}
\min \frac{Nf}{A}+\sum\limits_{i=1}^{N} {\frac{\kappa f A_i\sqrt{A_i}}{Ag(n_i)}}+\mu\left(A-\sum\limits_{i=1}^{N} {A_i}\right).\label{total_cost_system_Lag}
\end{equation}

By the first order condition, we have
\begin{equation*}
\mu= \frac{3}{2}\kappa f\left(A\sum\limits_{i=1}^{N} {g^2(n_i)}\right)^{-\frac{1}{2}}\text{ and }A_i=g^2(n_i)A\left(\sum\limits_{i=1}^{N} {g^2(n_i)}\right)^{-1},\forall i.
\end{equation*}

Substitute $A_i=\frac{g^2(n_i)A}{\sum_{i=1}^{N} {g^2(n_i)}},\forall i$ into~\eqref{total_cost_system} and we can get
\begin{align}
 z(N,A) \geq \frac{Nf}{A}+\sum\limits_{i=1}^{N} {\frac{\kappa f A^{\frac{1}{2}}g^2(n_i)}{\left(\sum\limits_{i=1}^{N} {g^2(n_i)}\right)^{\frac{3}{2}}}}=\frac{Nf}{A}+\frac{\kappa f A^{\frac{1}{2}}}{\left(\sum\limits_{i=1}^{N} {g^2(n_i)}\right)^{\frac{1}{2}}}
 \geq \frac{Nf}{A}+\frac{\kappa f A^{\frac{1}{2}}}{\sqrt{N}g\left(\sum_{i=1}^N n_i/N\right)}. \label{total_cost_system_opt}
\end{align}

The last inequality holds from \eqref{opt_lb}. \end{proof}

\subsection{Proof of Lemma~\ref{C_circumscribe_inbound_root}}\label{proof_C_circumscribe_inbound_root}
\begin{proof}
Since $\sin\alpha>0,\sin {\frac {\pi-2\alpha}{n_i-2}}>0$, $\forall\alpha\in\left[\frac{\pi}{n_i},\frac{\pi}{2}\right)$, we define a new function $\bar{H}(n_i,r,\alpha)=\frac{H(n_i,r,\alpha)}{\sin\alpha\sin {\frac {\pi-2\alpha}{n_i-2}}}$ such that $\bar{H}(n_i,r,\alpha)$ has the same root as $H(n_i,r,\alpha)$ with respect to $\alpha \in\left[\frac{\pi}{n_i}, \frac{\pi}{2}\right)$. From~\eqref{equation_alpha_1_2_V2_H}, we can express $\bar{H}(n_i,r,\alpha)$ as follows:
\begin{equation}
 \bar{H}(n_i,r,\alpha)=p\left(\frac {\pi-2\alpha}{n_i-2}\right)-p\left(\alpha\right)+q\left(\alpha\right),\frac{\pi}{n_i}\leq\alpha<\frac{\pi}{2},r \geq{0},n_i\geq 3,\label{equation_alpha_1_2_V3}
\end{equation}
where
\begin{align}
p(x)&=\frac{\cos^2x}{\sin x}\log(\tan\left(\frac{\pi}{4}+\frac{x}{2}\right)),0<x<\frac{\pi}{2},\label{define_p}\\
q(x)&=-2r \cos x-\frac{r(n_i-2)\cos {\frac {\pi-2x}{n_i-2}} \sin {\frac {\pi-2x}{n_i-2}} }{ \sin x},\frac{\pi}{n_i}\leq x<\frac{\pi}{2},r \geq{0},n_i\geq 3.\label{define_1}
\end{align}

It is easy to show that $\frac{dp(x)}{dx}<0$ when $0<x<\frac{\pi}{2}$. Thus $p(x)$ decreases strictly monotonically in the open interval $(0,\frac{\pi}{2})$, and the first two terms in \eqref{equation_alpha_1_2_V3} increases monotonically with $\alpha$.
By the same token, we find that $\frac{dq(x)}{dx}>0$ when $\frac{\pi}{n_i}\leq x< \frac{\pi}{2}$. Hence, $q(x)$ increases monotonically over $x$ when $\frac{\pi}{n_i}\leq x< \frac{\pi}{2}$.

As such, $\bar{H}(n_i,r,\alpha)$ increases monotonically over $\alpha \in\left[\frac{\pi}{n_i}, \frac{\pi}{2}\right)$. Meanwhile, as $H(n_i,r,\frac{\pi}{n_i})=-mn\cos \frac{\pi}{n_i}<0$ and $\bar{H}\left(n_i,r,\frac{\pi}{2}\right)=\frac{1}{2}$, the implicit equation $\bar{H}(n_i,r,\alpha)=0$ has one and only one root in the interval $\left[\frac{\pi}{n_i},\frac{\pi}{2}\right)$ for all $n_i\geq 3$ and $r \geq 0$. \end{proof}

\subsection{Proof of Lemma~\ref{P_circumscribe_inbound}}\label{proof_P_circumscribe_inbound}
\begin{proof}
We take the first order derivative of \eqref{z_all_inbound_outbound} with respect to $\alpha_i$, and the first order condition yields
\begin{align}
&(n_i-2)R_i^2F(\bar{\alpha}_i)\frac{d R_i}{d\alpha_i}+\frac{(n_i-2)R_i^3}{3}\frac{d F(\bar{\alpha}_i)}{d\bar{\alpha}_i}\frac{d \bar{\alpha}_i}{d\alpha_i}+2R_i^2F(\alpha_i)\frac{d R_i}{d\alpha_i}+\frac{2R_i^3}{3}\frac{d F(\alpha_i)}{d\alpha_i}+2\kappa rf A_i\cos\alpha_i\frac{d R_i}{d\alpha_i}\notag\\
&-2\kappa r f A_iR_i\sin\alpha_i=0,\label{z_all_inbound_outbound_Lag_De_R_2}
\end{align}
where $F(x)=\kappa f\cos^3x\int_{-x}^{x}\frac{1}{\cos^3 t} dt$ is defined in Appendix~\ref{proof_P_Given_Area_All}.

Since $(n_i-2)\bar{\alpha}_i+2\alpha_i=\pi$ always holds, then we have
\begin{equation}
\frac{d\bar{\alpha}_i}{d\alpha_i}=-\frac{2}{n_i-2},\label{inbound_outbound_de_alpha_1}
\end{equation}
and \eqref{z_all_inbound_outbound_Lag_De_R_2} can be simplified into the following
\begin{equation}
F(\alpha_i)-\frac{F(\bar{\alpha}_i)\sin\alpha_i\cos\alpha_i}{\sin\bar{\alpha}_i\cos\bar{\alpha}_i}+\kappa r f \cos\alpha_i \left[(n_i-2)\sin\bar{\alpha}_i\cos\bar{\alpha}_i+2\sin\alpha_i\cos\alpha_i\right]=0.\label{equation_alpha_1_2_V1}
\end{equation}

As $F(x)=\kappa f\cos^3x\int_{-x}^{x}\frac{1}{\cos^3 t} dt=\kappa f\cos^3x\left(\log\tan\left(\frac{x}{2}+\frac{\pi}{4}\right)+\frac{\tan x}{\cos x}\right)$, and let $\alpha=\alpha_i$, the left side of ~\eqref{equation_alpha_1_2_V1} can be rewritten as the following function
\begin{align}
H(n_i,r,\alpha)&= \sin\alpha\cos^2 {\frac {\pi-2\alpha}{n_i-2}} \log\tan \left( \frac{\pi}{4}+{\frac {\pi-2\alpha}{2(n_i-2)}} \right)-\cos^2\alpha\sin {\frac {\pi-2\alpha}{n_i-2}}\log \tan \left( \frac{\pi}{4}+\frac{\alpha}{2} \right)\notag \\
 &-r \sin 2\alpha\sin {\frac {\pi-2\alpha}{n_i-2}} -r (n_i-2)\cos {\frac {\pi-2\alpha}{n_i-2}} \sin^2 {\frac {\pi-2\alpha}{n_i-2}}.\label{equation_alpha_1_2_V2}
\end{align}

Simple algebra will show that the lower bound of \eqref{z_all_inbound_outbound} is
\allowdisplaybreaks\begin{align}
\frac{\kappa f A_i^{\frac{3}{2}}\left(1+\cos {\frac {\pi-2\alpha^*}{n_i-2}}\cot {\frac {\pi-2\alpha^*}{n_i-2}} \log\tan \left( \frac{\pi}{4}+{\frac {\pi-2\alpha^*}{2(n_i-2)}} \right)+4r \cos\alpha^*\right)}{3\left(2\cos \alpha^* \sin \alpha^* +(n_i-2)\cos {\frac {\pi-2\alpha^*}{n_i-2}} \sin {\frac {\pi-2\alpha^*}{n_i-2}}\right)^{\frac{1}{2}}},\label{z_all_inbound_outbound_1}
\end{align}
where $\alpha^*$ is the root of $H(n_i,r,\alpha)=0$. In light of Lemma~\ref{C_circumscribe_inbound_root}, this completes the proof. \end{proof}

\subsection{Proof of Proposition~\ref{Prop_lB_Z_N_A_inbound}}\label{proof_Prop_lB_Z_N_A_inbound}
\begin{proof}
We minimize the right hand side of~\eqref{total_cost_system_inbound} as an unconstrained optimization problem by adding constraint~\eqref{Area_con_system} into the objective with Lagrangian multiplier $\eta$; i.e.,
\begin{equation}
\frac{Nf}{A}+\sum\limits_{i=1}^{N} {\frac{\kappa f A_i\sqrt{A_i}}{3Ag(n_i,r )}}+\eta\left(A-\sum\limits_{i=1}^{N} {A_i}\right). \label{total_cost_system_inbound_Lag}
\end{equation}

By the first order conditions, we have
\begin{equation*}
\eta= \frac{3\kappa f}{2\sqrt{A\sum\limits_{i=1}^{N} {g^2(n_i,r )}}} \text{ and } A_i=\frac{g^2(n_i,r )A}{\sum\limits_{i=1}^{N} {g^2(n_i,r )}},\forall i.\label{eta_Lag}
\end{equation*}

Substitute $A_i=\frac{g^2(n_i,r )A}{\sum\limits_{i=1}^{N} {g^2(n_i,r )}},\forall i$ into~\eqref{total_cost_system_inbound}, we have
\begin{equation}
 z_{ub}(N,A) \geq \frac{Nf}{A}+\sum\limits_{i=1}^{N} {\frac{\kappa f A^{\frac{1}{2}}g^2(n_i,r )}{\left(\sum\limits_{i=1}^{N} {g^2(n_i,r )}\right)^{\frac{3}{2}}}}= \frac{Nf}{A}+\frac{\kappa f A^{\frac{1}{2}}}{\left(\sum\limits_{i=1}^{N} {g^2(n_i,r )}\right)^{\frac{1}{2}}},\forall \{n_i\},N,r. \label{total_cost_system_inbound_opt}
\end{equation}

Numerical examination of the first order and second order derivatives of smooth function $g^2(n,r )$ show that $\frac{\partial g^2(n,r )}{\partial n}>0,\frac{\partial^2 g^2(n,r )}{\partial n^2}<0$. Thus, $g^2(n,r )$ is concave and monotonically increasing over $n\in\left[3,+\infty\right)$. Thus,
\begin{equation}
 \sum\limits_{i=1}^{N} {g^2(n_i,r )}\leq{Ng^2\left(\sum\limits_{i=1}^{N}n_i/N,r \right)},,\forall \{n_i\},N,r.\label{opt_inbound_lb}
\end{equation}

Note that~\cite{Newman82} proved that $\sum_{i=1}^{N}n_i\leq{6N}$ for any usual tessellation of the plane. Thus, we further have
\begin{align}
\sum_{i=1}^{N} {\frac{A_i^{\frac{3}{2}}}{Ag(n_i,r )}}\geq \frac{A^{\frac{1}{2}}}{\sqrt{N}g\left(\sum_{i=1}^N n_i/N,r \right)}\geq \frac{A^{\frac{1}{2}}}{\sqrt{N}g\left(6,r \right)},\forall \{n_i\},N,r.\label{opt_inbound_lb1}
\end{align}
\end{proof}

\end{appendices}

\end{document}